\documentclass[12pt]{article}

\usepackage{amssymb}
\usepackage{amsmath}
\usepackage{cite}
\usepackage{microtype}
\usepackage{mathrsfs}

\begin{document}

\renewcommand{\citeleft}{{\rm [}}
\renewcommand{\citeright}{{\rm ]}}
\renewcommand{\citepunct}{{\rm,\ }}
\renewcommand{\citemid}{{\rm,\ }}

\newcommand{\LegendreScalarProduct}[3][n]{\left[ #2, #3 \right]_{#1}}

\newcounter{abschnitt}
\newtheorem{satz}{Theorem}
\newtheorem{theorem}{Theorem}[abschnitt]
\newtheorem{koro}[theorem]{Corollary}
\newtheorem{prop}[theorem]{Proposition}
\newtheorem{lem}[theorem]{Lemma}
\newtheorem{approp}{Proposition}[section]
\newtheorem{apthm}[approp]{Theorem}

\newcounter{saveeqn}
\newcommand{\alpheqn}{\setcounter{saveeqn}{\value{abschnitt}}
\renewcommand{\theequation}{\mbox{\arabic{saveeqn}.\arabic{equation}}}}
\newcommand{\reseteqn}{\setcounter{equation}{0}
\renewcommand{\theequation}{\arabic{equation}}}

\hyphenpenalty=7500

\sloppy

\phantom{a}

\vspace{-2.2cm}

\begin{center}
\begin{Large} {Log-Concavity Properties of Minkowski Valuations} \\[0.7cm] \end{Large}

\begin{large} Astrid Berg, Lukas Parapatits, \\[0.2cm] Franz E.\
Schuster, Manuel Weberndorfer \end{large}
\end{center}

\vspace{-1.1cm}

\begin{quote}
\footnotesize{ \vskip 1truecm\noindent {\bf Abstract.} New Orlicz
Brunn--Minkowski inequalities are established for rigid motion
compatible Minkowski valuations of arbitrary degree. These extend
classical log-concavity properties of intrinsic volumes and
generalize seminal results of Lutwak and others. Two different
approaches which refine previously employed techniques are
explored. It is shown that both lead to the same class of
Minkowski valuations for which these inequalities hold. An
appendix by Semyon Alesker contains the proof of a new classification of
generalized translation invariant valuations.}
\end{quote}

\vspace{0.6cm}

\centerline{\large{\bf{ \setcounter{abschnitt}{1}
\arabic{abschnitt}. Introduction}}} \alpheqn

\vspace{0.6cm}

The fundamental log-concavity property of the volume functional
is \linebreak expressed by the multiplicative form of the
Brunn--Minkowski inequality:
\begin{equation} \label{logbmineq}
V_n((1-\lambda) K + \lambda L) \geq
V_n(K)^{1-\lambda}V_n(L)^{\lambda},
\end{equation}
where $K$ and $L$ are convex bodies (non-empty compact convex
sets) in $\mathbb{R}^n$ with non-empty interiors, $0 < \lambda <
1$, and $+$ denotes Minkowski addition. Equality holds in
(\ref{logbmineq}) if and only if $K$ and $L$ are translates of
each other. The excellent survey of Gardner
\textbf{\cite{gardner02}} gives a comprehensive overview of
different aspects and consequences of the Brunn--Minkowski
inequality. Here we just mention that it directly implies the
classical {\it Euclidean} isoperimetric inequality.

Projection bodies of convex bodies were defined at the turn of the
previous century by Minkowski. In 1984 Lutwak
\textbf{\cite{lutwak84}} discovered that an {\it affine}
isoperimetric inequality of Petty \textbf{\cite{petty67}} for
{\it polar} projection bodies is not only significantly stronger
than the Euclidean isoperimetric inequality, but in fact an
optimal version of this classical inequality. For the tremendous
impact of Petty's inequality and its generalizations see,
e.g., \textbf{\cite{habschu09, LYZ2000a, LYZ2010, tuo12,
Zhang99}}. The problem of finding sharp bounds for the volume of
projection bodies, given the volume of the original body, remains
a central quest in convex geometric analysis. It has led, among
many other results, to the discovery of important log-concavity
properties of the volume of projection bodies. In fact, Lutwak
\textbf{\cite{lutwak93}} established not only Brunn--Minkowski
type inequalities for the volume of projection bodies, but for
all the {\it intrinsic volumes} of projection bodies of arbitrary
order (see Section 2).

\pagebreak

In the present article we investigate a common generalization of
Lutwak's Brunn--Minkowski inequalities for projection bodies and
inequality (\ref{logbmineq}), more specifically, its version for
all the intrinsic volumes. To be more precise, \linebreak we
establish new log-concavity properties of intrinsic volumes of
convex body valued \emph{valuations} which intertwine rigid motions.
This line of research has its origins in the discovery of the
special place of projection bodies in affine \linebreak geometry:
Ludwig \textbf{\cite{ludwig02, Ludwig:Minkowski}} characterized
the projection body operator as the unique continuous Minkowski
valuation which is translation invariant and $\mathrm{GL}(n)$
contravariant (see \textbf{\cite{haberl12, Ludwig:matrix,
Ludwig10a, parapTAMS13, parapJLMS13, SchuWan11, wannerer10}} for
related results).

In recent years, it has become apparent that several geometric
inequalities for projection bodies and, more general, valuations
intertwining the group of affine transformations, in fact, hold
for much larger classes of valuations intertwining merely rigid
motions. First such results were obtained in
\textbf{\cite{Schu06b}}, where the Brunn--Minkowski inequalities
for projection bodies of Lutwak were generalized to translation
invariant and $\mathrm{SO}(n)$ \emph{equivariant} Minkowski
valuations of degree $n - 1$. Although considerable efforts have
been invested ever since to show that these log-concavity
properties extend to Minkowski valuations of arbitrary degree (see
\textbf{\cite{ABS2011, papschu12, Schu09}}), the conjectured
complete family of inequalities has only partially been obtained
(compare Section 2).

For the inequalities established so far two different approaches
were used. While in \textbf{\cite{Schu06b}} and
\textbf{\cite{Schu09}} integral representations of (even)
Minkowski valuations \linebreak which are translation invariant
and $\mathrm{SO}(n)$ equivariant were crucial, in
\textbf{\cite{papschu12}} the Hard Lefschetz \emph{derivation} operator
on Minkowski valuations \textbf{\cite{Alesker03, papschu12}},
\linebreak together with a symmetry property of bivaluations
\textbf{\cite{ABS2011}}, was the key ingredient. \linebreak In
this paper we show that the Hard Lefschetz {\it integration}
operator on Minkowski valuations \textbf{\cite{Alesker04,
Bernig07, SchuWan13}} on one hand and a recent representation
theorem for Minkowski valuations \textbf{\cite{SchuWan13,
SchuWan14, wannerer12}} on the other hand lead to a natural class
of Minkowski valuations which exhibit log-concavity properties.
All the Brunn--Minkowski inequalities for Minkowski valuations
established before turn out to be special cases of our new
results. From new monotonicity \linebreak properties of these
Minkowski valuations, we are  able to deduce a complete
characterization of equality cases without any smoothness
assumptions that were required before.

Moreover, all previously obtained and new Brunn--Minkowski
inequalities for Minkowski valuations are shown to not only hold
for Minkowski addition but for all \emph{commutative} Orlicz Minkowski
additions (introduced in \textbf{\cite{gardhugweil13}}) of convex bodies. This
includes, in particular, all $L_p$ Minkowski additions.

\pagebreak

\centerline{\large{\bf{ \setcounter{abschnitt}{2}
\arabic{abschnitt}. Statement of principal results}}}

\reseteqn \alpheqn

\vspace{0.6cm}

Let $\mathcal{K}^n$ denote the space of convex bodies in
$n$-dimensional Euclidean space $\mathbb{R}^n$ endowed with the
Hausdorff metric. Throughout the article we assume that $n \geq
3$. A convex body $K$ is uniquely determined by its support
function $h(K,u) = \max\{u \cdot x: x \in K\}$ for $u \in
S^{n-1}$. For $i \in \{0, \ldots, n\}$, let $V_i(K)$ denote the
$i$th intrinsic volume of $K$ (see Section 3).

A map $\Phi: \mathcal{K}^n \to \mathcal{K}^n$ is called a {\it
Minkowski valuation} if
\[\Phi K + \Phi L = \Phi(K \cup L) + \Phi(K \cap L)  \]
whenever $K \cup L \in \mathcal{K}^n$ and addition on
$\mathcal{K}^n$ is Minkowski addition.

The theory of {\it scalar} valued valuations has long played a
prominent role in convex geometry (see, e.g.,
\textbf{\cite{hadwiger51, Klain:Rota}} for the history of scalar
valuations and \linebreak \textbf{\cite{Alesker01,bernigfu10, fu06, habparap13, centro, papwann12, wannerer13}} for
more recent results). Systematic investigations of Minkowski
valuations have only been initiated about a decade ago by Ludwig
\textbf{\cite{ludwig02, Ludwig:matrix, Ludwig:Minkowski}}. These valuations
arise naturally from data about projections and sections of
convex bodies and form an integral part of geometric tomography.
As first examples we mention here the projection body maps $\Pi_i:
\mathcal{K}^n \to \mathcal{K}^n$ of order $i \in \{1, \ldots, n -
1\}$, defined by
\[h(\Pi_iK,u) = V_i(K|u^\bot), \qquad u \in S^{n-1}.  \]

While the entire family $\Pi_i$ is translation invariant and
$\mathrm{SO}(n)$ equivariant, \linebreak the classic projection
body map $\Pi_{n-1}$ is the only one among them which
inter\-twines {\it linear} transformations (see
\textbf{\cite{ludwig02}}). In fact, there is only a small number
of Minkowski valuations which are compatible with affine
transformations (see \textbf{\cite{abardia12, abardiabernig, haberl12, Ludwig:Minkowski,
SchuWan11, wannerer10}} for their classification).

In this article we establish new log-concavity properties for the
class $\mathbf{MVal}_j$ of continuous, translation invariant {\it
and} $\mathrm{SO}(n)$ equivariant Minkowski valuations of a given
degree $j$ of homogeneity (by a result of McMullen
\textbf{\cite{McMullen77}}, only integer degrees $0 \leq j \leq
n$ can occur; cf.\ Section 5). A first such result was obtained
by Lutwak \textbf{\cite[\textnormal{Theorem 6.2}]{lutwak93}} for
projection bodies of arbitrary order. In an equivalent multiplicative form it states the following: If $K, L \in \mathcal{K}^n$
have non-empty interiors, $1 \leq i \leq n$, and $2 \leq j \leq n
- 1$, then for all $\lambda \in (0,1)$,
\begin{equation} \label{piiinequ}
V_i(\Pi_j((1-\lambda)K + \lambda L)) \geq
V_i(\Pi_jK)^{1-\lambda}V_i(\Pi_jL)^{\lambda},
\end{equation}
with equality if and only if $K$ and $L$ are translates of each
other.

Inequalities (\ref{piiinequ}) have been generalized in different
directions: Abardia and Bernig \textbf{\cite{abardiabernig}}
extended (\ref{piiinequ}) to the entire class of {\it complex}
projection bodies. Analogues of (\ref{piiinequ}) were established
in \textbf{\cite{Schu06b}} for all valuations in
$\mathbf{MVal}_{n-1}$ and then in \textbf{\cite{Schu09}} for {\it
even} valuations in $\mathbf{MVal}_j$ in the case $i = j + 1$.
The assumption on the parity could later be omitted in
\textbf{\cite{ABS2011}}. In the Euclidean setting, the most general
result to date can be stated (in multiplicative form) as follows, where
we call the Minkowski valuation which maps every convex body to
the set containing only the origin {\it trivial}.

\begin{satz} \label{thm1} {\bf (\!\!\cite{papschu12})} Let $\Phi_j \in \mathbf{MVal}_j$, $2 \leq j \leq n - 1$, be non-trivial. If $K, L \in \mathcal{K}^n$ and $1
\leq i \leq j + 1$, then for all $\lambda \in (0,1)$,
\begin{equation} \label{desinequ}
V_i(\Phi_j((1-\lambda)K+\lambda L)) \geq
V_i(\Phi_jK)^{1-\lambda}V_i(\Phi_jL)^\lambda.
\end{equation}
If $K$ and $L$ are of class $C^2_+$, then equality holds if and
only if $K$ and $L$ are translates of each other.
\end{satz}

Note that Theorem \ref{thm1} establishes (\ref{desinequ}) only
for $1 \leq i \leq j + 1$, while in Lutwak's family of
inequalities (\ref{piiinequ}) the range of $i$ does not depend on
$j$. The proof of Theorem \ref{thm1} used ideas from
\textbf{\cite{ABS2011}} and the existence of a new derivation
operator $\Lambda$ on Minkowski valuations established in
\textbf{\cite{papschu12}} (see also Section 5). For $\Phi \in
\mathbf{MVal}_j$, there exists $\Lambda \Phi \in
\mathbf{MVal}_{j-1}$ such that
\[ h((\Lambda \Phi)(K),\cdot) = \left . \frac{d}{dt} \right |_{t=0} h(\Phi(K + tB),\cdot).\]
This definition was motivated by a similar derivation operator
introduced by Alesker \textbf{\cite{Alesker03}} in the theory of {\it scalar
valued} valuations. There it is widely used to deduce results for
valuations of degree $i$ from those for valuations of \linebreak
some degree $j
> i$. The key to the proof of Theorem~\ref{thm1} was the
following generalization of a symmetry property of bivaluations
obtained in \textbf{\cite{ABS2011}}.

\begin{satz} \label{thm2} {\bf (\!\!\cite{papschu12})} Let $\Phi_j \in \mathbf{MVal}_j$, $2 \leq j \leq n - 1$. If $1 \leq i \leq j + 1$, then
\begin{equation} \label{durch1}
W_{n-i}(K,\Phi_jL)=\frac{(i-1)!}{j!}\,W_{n-j-1}(L,(\Lambda^{j+1-i}\Phi_j)(K))
\end{equation}
for every $K, L \in \mathcal{K}^n$.
\end{satz}

Here $W_{m}(K,L)$ denotes the mixed volume $V(K[n-m-1],B[m],L)$
with $n-m-1$ copies of $K$ and $m$ copies of the Euclidean unit
ball $B$. We will see in Section 6 that Theorem \ref{thm2}
follows directly from a recently obtained integral
representation of Minkowski valuations intertwining rigid motions.

\pagebreak

An obvious idea for a proof of (\ref{desinequ}) for the remaining
cases $j + 2 \leq i \leq n$ is to establish a counterpart of
Theorem \ref{thm2} for the Hard Lefschetz integration operator:
For $\Phi \in \mathbf{MVal}_j$, there exists $\mathfrak{L} \Phi
\in \mathbf{MVal}_{j+1}$ such that
\[ h((\mathfrak{L} \Phi)(K),\cdot) = \int_{\mathrm{AGr}_{n-1,n}} h(\Phi(K \cap E),\cdot)\,dE,  \]
where $\mathrm{AGr}_{n-1,n}$ denotes the affine Grassmannian of
$n - 1$ planes in $\mathbb{R}^n$ and where we integrate with
respect to the suitably normalized invariant measure on
$\mathrm{AGr}_{n-1,n}$ (see Section 5). For {\it scalar valued}
valuations the operator $\mathfrak{L}$ was first defined in
\textbf{\cite{Alesker04}} and used to deduce results for
valuations of degree~$i$ from those for valuations of some degree
$j < i$. As an operator on Minkowski valuations, $\mathfrak{L}$
was first considered in \textbf{\cite{SchuWan13}}.

Our first result is a version of Theorem \ref{thm2} for the
operator $\mathfrak{L}$. However, the situation is more delicate
in this case and we will see that a full analogue of
(\ref{durch1}) only holds for a subclass of Minkowski valuations.
To define this class let $\mathbf{MVal}_j^{\infty}$ denote the
set of translation invariant and $\mathrm{SO}(n)$ equivariant
{\it smooth} Minkowski valuations (cf.\ Section 5).

\vspace{0.4cm}

\noindent {\bf Definition}  \emph{For $1 \leq i,j \leq n - 1$, let
$\mathbf{MVal}^{\infty}_{j,i} \subseteq \mathbf{MVal}_j^\infty$ be
defined by
\[\mathbf{MVal}^{\infty}_{j,i} = \left \{ \begin{array}{ll} \Lambda^{i-j}(\mathbf{MVal}_i^\infty) & \mbox{ if } i > j, \\ \mathbf{MVal}_j^\infty & \mbox{ if } i \leq j. \end{array} \right .  \]
We write $\mathbf{MVal}_{j,i}$ for the closure of
$\mathbf{MVal}_{j,i}^{\infty}$ in the topology of uniform
convergence on compact subsets.}

\vspace{0.3cm}

We will see in Section 5 that $\Lambda: \mathbf{MVal}_j^{\infty}
\rightarrow \mathbf{MVal}_{j-1}^{\infty}$ is injective for
\linebreak $2 \leq j \leq n$. Thus, for $i
> j$, the inverse map $(\Lambda^{i-j})^{-1}: \mathbf{MVal}^{\infty}_{j,i} \rightarrow \mathbf{MVal}_i^{\infty}$ is well defined and will be denoted by $\Lambda^{j-i}$.

Our counterpart of Theorem \ref{thm2} can now be stated as
follows:

\begin{satz} \label{thm3} Let $\Phi_j \in
\mathbf{MVal}_j^{\infty}$, $2 \leq j \leq n - 1$. For $j + 2 \leq
i \leq n$ and every convex body $L \in \mathcal{K}^n$, there
exists a \emph{generalized valuation} $\gamma_{i,j}(L,\cdot) \in
\mathbf{Val}_1^{-\infty}$ such that
\begin{equation*}
W_{n-i}(K,\Phi_jL)=\gamma_{i,j}(L,(\mathfrak{L}^{i-j-1}\Phi_j)(K))
\end{equation*}
for every $K \in \mathcal{K}^n$. Moreover, if $\Phi_j \in
\mathbf{MVal}_{j,i-1}^{\infty}$, then
\[\gamma_{i,j}(L,(\mathfrak{L}^{i-j-1}\Phi_j)(K)) = \frac{(i-1)!}{j!} W_{n-1-j}(L,(\Lambda^{j+1-i}\Phi_j)(K)).   \]
\end{satz}

Generalized translation invariant valuations were introduced recently by Alesker and Faifman \textbf{\cite{AleskerFaifman}}. We recall their definition and basic properties (in
particular, of the space $\mathbf{Val}_1^{-\infty}$) in Section 5. A crucial ingredient in the proof of Theorem \ref{thm3}
is a new classification of generalized valuations from $\mathbf{Val}_1^{-\infty}$. We are very grateful to Semyon Alesker for
communicating to us a proof of this result and his permission to include it in an appendix of this article.

Using Theorem \ref{thm3}, we establish in Section 6 the main
result of this article:

\begin{satz} \label{thm4} Let $1 \leq i \leq n$ and let $\Phi_j \in \mathbf{MVal}_{j,i-1}$, $2 \leq j \leq n - 1$, be non-trivial. If $K, L \in \mathcal{K}^n$ have non-empty interiors,
then for all $\lambda \in (0,1)$,
\begin{equation} \label{thm4form}
V_i(\Phi_j((1-\lambda)K+\lambda L)) \geq
V_i(\Phi_jK)^{1-\lambda}V_i(\Phi_jL)^\lambda,
\end{equation}
with equality if and only if $K$ and $L$ are translates of each
other.
\end{satz}

Since $\mathbf{MVal}_{j,i-1} = \mathbf{MVal}_j$ for $i \leq j +
1$, Theorem \ref{thm4} includes both \linebreak Lutwak's
inequalities (\ref{piiinequ}) and Theorem~\ref{thm1} as special
cases. Also note that the smoothness assumption for the bodies
$K$ and $L$ in the equality conditions of (\ref{desinequ}) is no
longer required. This follows from new monotonicity properties of
the Minkowski valuations in $\mathbf{MVal}_{j,i-1}$, which we
prove in Section 6.

In the last part of the article we explain that our proof of
Theorem \ref{thm4} can be modified to yield an even stronger
result. More precisely, we show that (\ref{thm4form}) not only holds for
the usual Minkowski addition but, in fact, for all
commutative Orlicz Minkowski additions introduced by Gardner, Hug,
and Weil \textbf{\cite{gardhugweil13}}. In particular, this
includes all the $L_p$ Minkowski additions.

Let $\Theta_1$ denote the set of convex functions $\varphi:
[0,\infty) \rightarrow [0,\infty)$ satisfying $\varphi(0)=0$ and
$\varphi(1)=1$. For $\varphi \in \Theta_1$ and $K, L \in
\mathcal{K}^n$ containing the origin, we write $K
+_{\varphi,\lambda} L$ for the Orlicz Minkowski convex combination of $K$
and $L$ (see Section 3 for the definition).

\begin{satz} \label{mainorlicz} Let $\varphi \in \Theta_1$, $1 \leq i \leq n$, and let $\Phi_j \in
\mathbf{MVal}_{j,i-1}$, $2 \leq j \leq n - 1$, be non-trivial. If
$K, L \in \mathcal{K}^n$ contain the origin, then for all
$\lambda \in (0,1)$,
\begin{equation} \label{yoda1}
V_i(\Phi_j(K+_{\varphi,\lambda} L)) \geq
V_i(\Phi_jK)^{1-\lambda}V_i(\Phi_jL)^\lambda.
\end{equation}
When $\varphi$ is strictly convex and $K$ and $L$ have non-empty
interiors, equality holds if and only if $K = L$.
\end{satz}

We will explain in Section 3 that by a recent result of Gardner,
Hug, and Weil \textbf{\cite{gardhugweil13}} (Theorem
\ref{charcomorlicz} below), inequality (\ref{yoda1}) holds for all
\emph{commutative} Orlicz Minkowski additions.

\pagebreak

\centerline{\large{\bf{ \setcounter{abschnitt}{3}
\arabic{abschnitt}. Background material on convex bodies}}}

\reseteqn \alpheqn \setcounter{theorem}{0}

\vspace{0.6cm}

For quick later reference we collect in this section some basic facts from convex geometry, in particular, on additions of convex bodies and inequalities for mixed volumes. As
general reference for this material we recommend the
book by Schneider \textbf{\cite{schneider93}} and the article \textbf{\cite{gardhugweil13}}.

For a convex body $K \in \mathcal{K}^n$, the definition of the support function $h(K,u)=\max\{u\cdot x: x
\in K\}$, $u \in S^{n-1}$, implies $h(\vartheta K,u)=h(K, \vartheta^{-1}u)$ for every $u
\in S^{n-1}$ and $\vartheta \in \mathrm{SO}(n)$. Since every twice continuously differentiable function on
$S^{n-1}$ is a difference of support functions (see, e.g.,
\textbf{\cite[\textnormal{p.\ 49}]{schneider93}}), the subspace
spanned by support functions $\{h(K,\cdot) - h(L,\cdot): K, L \in
\mathcal{K}^n\}$ is dense in $C(S^{n-1})$. The {\it Steiner point} $s(K)$ of $K \in \mathcal{K}^n$ is defined by
\[s(K)=\frac{1}{\kappa_n}\int_{S^{n-1}}h(K,u)u\,du.  \]
Here and in the following we use $du$ to denote integration with respect to spherical Lebesgue measure and $\kappa_m$ for
the $m$-dimensional volume of the unit ball in $\mathbb{R}^m$. The Steiner point map is the
unique vector valued, rigid motion equivariant and continuous
valuation on $\mathcal{K}^n$ (see e.g.,
\textbf{\cite[\textnormal{p.\ 363}]{schneider93}}).

For $K, L \in \mathcal{K}^n$ and $s, t \geq 0$, the support
function of the Minkowski combination $s\,K + t\,L$ is given by
\[h(s\,K + t\,L,\cdot)=s\,h(K,\cdot)+t\,h(L,\cdot).\]

On the set $\mathcal{K}_{\mathrm{o}}^n$ of convex bodies
containing the origin, Firey introduced in the 1960s a more
general way of combining convex sets. For $K, L \in
\mathcal{K}^n_{\mathrm{o}}$,  $s,t \geq 0$, and $1 \leq p <
\infty$, the \emph{$L_p$ Minkowski combination} $s \cdot K +_p t
\cdot L$ is defined by
\[h(s \cdot K +_p t \cdot L,\cdot)^p = s\, h(K,\cdot)^p + t\, h(L,\cdot)^p.  \]

Initiated by Lutwak \textbf{\cite{Lutwak:1993, Lutwak:1996}}, in the last two decades an entire $L_p$ theory of convex bodies was developed which
represents a powerful extension of the classical Brunn--Minkowski theory (see, e.g.,
\textbf{\cite{habschu09, LYZ2000a, LYZ200b, parapJLMS13, parapTAMS13, SchuWeb12, webern13}}).

A still more recent extension of the Brunn--Minkowski theory goes
back to two articles of Lutwak, Yang, and Zhang
\textbf{\cite{LYZ2010, LYZ2010a}} and a paper by Haberl, Lutwak,
Yang, and Zhang \textbf{\cite{HLYZ}}. While these articles form
the starting point of an emerging Orlicz Brunn--Minkowski theory
that generalizes the $L_p$ theory of convex bodies in the same way
that Orlicz spaces generalize $L_p$ spaces, the
fundamental notion of an Orlicz Minkowski combination of convex
bodies was introduced later by Gardner, Hug, and Weil
\textbf{\cite{gardhugweil13}}.

\pagebreak

As before let $\Theta_1$ be the set of convex functions $\varphi:
[0,\infty) \rightarrow [0,\infty)$ satisfying $\varphi(0)=0$ and
$\varphi(1)=1$. For $K, L \in \mathcal{K}^n_{\mathrm{o}}$, $s, t
\geq 0$, and $\varphi, \psi \in \Theta_1$, the \emph{Orlicz
Minkowski combination} $+_{\varphi,\psi}(K,L,s,t)$ is defined by
\[h(+_{\varphi,\psi}(K,L,s,t),u) = \inf \left \{\alpha > 0: s\, \varphi \left (\frac{h(K,u)}{\alpha} \right ) + t\, \psi \left (\frac{h(L,u)}{\alpha} \right ) \leq 1  \right \} \]
for $u \in S^{n-1}$. The notation $+_{\varphi,\psi}(K,L,s,t)$ is
necessitated by the fact that it is not possible in general to
isolate an Orlicz scalar multiplication. We note that for
$\varphi(t) = \psi(t) = t^p$, $p \geq 1$, the Orlicz Minkowski
combination $+_{\varphi,\psi}(K,L,s,t)$ equals the $L_p$
Minkowski combination $s \cdot K +_p t \cdot L$.

For $s = t = 1$, we write $K +_{\varphi,\psi} L$ instead of
$+_{\varphi,\psi}(K,L,1,1)$ and call this the Orlicz Minkowski sum of
$K$ and $L$. In fact, Gardner, Hug, and Weil defined a more
general Orlicz addition but proved (see
\textbf{\cite[\textnormal{Theorem 5.5}]{gardhugweil13}}) that
their definition leads (essentially) to the Orlicz Minkowski
addition as defined here \emph{and} the $L_{\infty}$ Minkowski
addition obtained as the Hausdorff limit of the $L_p$~Minkowski
addition, that is, for $K, L \in \mathcal{K}^n_{\mathrm{o}}$,
\[K +_{\infty} L = \lim \limits_{p \rightarrow \infty} K +_p L = \mathrm{conv}(K \cup L).  \]

While all $L_p$ Minkowski additions are commutative, in general,
the Orlicz Minkowski addition of convex bodies is not. A classification of
those Orlicz additions which are commutative was obtained by
Gardner, Hug, and Weil.

\begin{theorem} {\bf ($\!\!$\cite{gardhugweil13})} \label{charcomorlicz}
Let $\varphi, \psi \in \Theta_1$. The addition $+_{\varphi,\psi}:
\mathcal{K}^n_{\mathrm{o}} \times \mathcal{K}^n_{\mathrm{o}}
\rightarrow \mathcal{K}^n_{\mathrm{o}}$ is commutative if and
only if there exists $\phi \in \Theta_1$ such that
$+_{\varphi,\psi} = +_{\phi,\phi}$.
\end{theorem}

In the following we will only be interested in commutative Orlicz
additions. For $K, L \in \mathcal{K}^n_{\mathrm{o}}$,
$\varphi \in \Theta_1$, and $\lambda \in (0,1)$ we use $K
+_{\varphi,\lambda} L$ to denote the Orlicz Minkowski convex
combination $+_{\varphi,\varphi}(K, L,(1-\lambda),\lambda)$. More
explicitly,
\begin{equation*}
h(K\! +_{\varphi,\lambda}\! L,u) = \inf\! \left \{\!\alpha > 0\!:
(1-\lambda) \varphi\! \left (\!\frac{h(K,u)}{\alpha}\! \right ) +
\lambda \varphi\! \left (\!\frac{h(L,u)}{\alpha}\! \right ) \leq
1\! \right \}
\end{equation*}
for $u \in S^{n-1}$. For the proof of Theorem \ref{mainorlicz} we
need the following simple fact.

\begin{lem} \label{orliczinclusion} If $\varphi \in \Theta_1$ and $K, L \in \mathcal{K}^n_{\mathrm{o}}$,
then for all $\lambda \in (0,1)$,
\begin{equation}\label{eq:incl}
K +_{\varphi,\lambda}  L \supseteq (1-\lambda)K + \lambda  L.
\end{equation}
\end{lem}
{\it Proof.} For $u \in S^{n-1}$ choose $t > h(K
+_{\varphi,\lambda} L,u)$. Then, by the convexity of $\varphi$ and
the definition of $K +_{\varphi,\lambda} L$, we have
\[\varphi\!\left(\frac{(1-\lambda)h(K,u)+\lambda h(L,u)}{t}\right)
\leq (1-\lambda)\varphi\!\left(\frac{h(K,u)}{t}\right) +
\lambda\varphi\!\left(\frac{h(L,u)}{t}\right) \leq 1.\] Since every $\varphi \in \Theta_1$ is
increasing and satisfies \mbox{$\varphi(1)=1$}, we conclude that
\[(1-\lambda)h(K,u)+\lambda h(L,u)\leq t.\]
Now, letting $t$ approach $h(K +_{\varphi,\lambda} L,u)$, we
obtain the desired inclusion (\ref{eq:incl}). $\phantom{aaa}$
\hfill $\blacksquare$

\vspace{0.3cm}

By a classical result of Minkowski, the volume
of a Minkowski linear combination $\lambda_1K_1 +
\cdots + \lambda_mK_m$, where $K_1, \ldots, K_m \in
\mathcal{K}^n$ and $\lambda_1, \ldots,
\lambda_m \geq 0$, can be expressed as a homogeneous polynomial of
degree $n$,
\begin{equation} \label{mixed}
V_n(\lambda_1K_1 + \cdots +\lambda_m K_m)=\sum
\limits_{j_1,\ldots, j_n=1}^m
V(K_{j_1},\ldots,K_{j_n})\lambda_{j_1}\cdots\lambda_{j_n},
\end{equation}
where the coefficients $V(K_{j_1},\ldots,K_{j_n})$, called {\it
mixed volumes} of $K_{j_1}, \ldots, K_{j_n}$, depend only on $K_{j_1}, \ldots, K_{j_n}$ and are symmetric in their
arguments. For $K, L \in \mathcal{K}^n$ and $0 \leq i \leq n$, we denote
the mixed volume with $i$ copies of $K$ and $n - i$ copies of $L$ by $V(K[i],L[n-i])$. For $K, K_1, \ldots, K_i \in
\mathcal{K}^n$ \linebreak and $\mathbf{C}=(K_1,\ldots,K_i)$, we write
$V_i(K,\mathbf{C})$ instead of
$V(K,\ldots,K,K_1,\ldots,K_i)$.

The mixed volume $W_i(K,K) =: W_i(K)$ is called the \emph{$i$th
quermassintegral} of $K$. The \emph{$i$th intrinsic volume}
$V_i(K)$ of $K$ is defined by
\begin{equation*} \label{viwi}
\kappa_{n-i}V_i(K)=\binom{n}{i} W_{n-i}(K).
\end{equation*}
A special case of (\ref{mixed}) is
the classical {\it Steiner formula} for the volume of the
parallel set of $K$ at distance $r > 0$,
\[V(K + rB)=\sum \limits_{i=0}^n r^i{n \choose i}W_i(K) = \sum
\limits_{i=0}^n r^{n-i}\kappa_{n-i}V_i(K).\]

A fundamental inequality for mixed volumes is the general
Minkowski inequality (see \textbf{\cite[\textnormal{p.\
427}]{schneider93}}): If $2 \leq i \leq n$ and $K, L \in
\mathcal{K}^n$ have dimension at least $i$, then
\begin{equation} \label{genmink}
W_{n-i}(K,L)^i \geq W_{n-i}(K)^{i-1}W_{n-i}(L),
\end{equation}
with equality if and only if $K$ and $L$ are homothetic.

\pagebreak

A consequence of (\ref{genmink}) and the homogeneity of
quermassintegrals is the (multiplicative) Brunn--Minkowski
inequality: If $2 \leq i \leq n$ and $K, L \in \mathcal{K}^n$
have dimension at least $i$, then for all $\lambda \in (0,1)$,
\begin{equation} \label{quermassbm}
W_i((1-\lambda)K+\lambda L) \geq
W_i(K)^{1-\lambda}W_i(L)^{\lambda},
\end{equation}
with equality if and only if $K$ and $L$ are translates of each
other.

A further generalization of inequality (\ref{quermassbm}) (where
the equality conditions are not yet known) is the following: If
$0 \leq i \leq n-2$, $K, L, K_1, \ldots, K_i \in \mathcal{K}^n$
and $\mathbf{C}=(K_1,...,K_i)$, then for all $\lambda \in (0,1)$,
\begin{equation} \label{mostgenbm}
V_i((1-\lambda)K + \lambda L,\mathbf{C}) \geq
V_i(K,\mathbf{C})^{1-\lambda}V_i(L,\mathbf{C})^{\lambda}.
\end{equation}

Associated with a convex body $K \in \mathcal{K}^n$ is a family
of Borel measures $S_i(K,\cdot)$, $0 \leq i \leq n - 1$, on
$S^{n-1}$, called the {\it area measures of order} $i$ of $K$.
They are uniquely determined by the property that
\begin{equation} \label{defsi}
W_{n-1-i}(K,L)=\frac{1}{n}\int_{S^{n-1}} h(L,u)\,dS_i(K,u)
\end{equation}
for all $L \in \mathcal{K}^n$. If $K \in \mathcal{K}^n$ has
non-empty interior, then, by a theorem of
Aleksandrov--Fenchel--Jessen (see, e.g.,
\textbf{\cite[\textnormal{p.\ 449}]{schneider93}}), each of the
measures $S_i(K,\cdot)$, $1 \leq i \leq n - 1$, determines $K$ up
to translations.

For $1 \leq j \leq n - 1$ and $r > 0$, we have the Steiner type
formula
\begin{equation*} \label{steinsi}
S_j(K + rB,\cdot)=\sum_{i=0}^j r^{j-i} {j \choose i}
S_{i}(K,\cdot).
\end{equation*}

A body $K \in \mathcal{K}^n$ is of class $C^2_+$ if the boundary
of $K$ is a $C^2$ submanifold of $\mathbb{R}^n$ with everywhere
positive curvature. In this case, each measure $S_i(K,\cdot)$, $0 \leq i \leq n - 1$, is
absolutely continuous with respect to spherical Lebesgue measure and its density is (up to a constant)
given by the $i$th elementary symmetric function of the principal
radii of curvature of $K$.

The center of mass (centroid) of every area measure of a convex
body is at the origin, that is, for every $K \in \mathcal{K}^n$
and all $i \in \{0, \ldots, n - 1\}$, we have
\begin{equation*}
\int_{S^{n-1}} u\,dS_i(K,u) = o.
\end{equation*}
The set $\mathcal{S}_i$ of all area measures of order $i$ of
convex bodies in $\mathcal{K}^n$ is dense in the set of all
non-negative finite Borel measures on $S^{n-1}$ with centroid at
the origin, endowed with the weak topology, if and only if $i =
n-1$. However, $\mathcal{S}_i - \mathcal{S}_i$, $1 \leq i \leq n - 1$, is dense in the
set $\mathcal{M}_{\mathrm{o}}(S^{n-1})$ of all \emph{signed}
finite Borel measures on $S^{n-1}$ with centroid at the origin
(see, e.g., \textbf{\cite[\textnormal{p.\ 477}]{schneider93}}).

\pagebreak

\centerline{\large{\bf{ \setcounter{abschnitt}{4}
\arabic{abschnitt}. Spherical harmonics and distributions}}}

\reseteqn \alpheqn \setcounter{theorem}{0}

\vspace{0.6cm}

In this section we collect facts about spherical harmonics, in particular,
on the series expansion of distributions on the sphere. We also recall C.\ Berg's \linebreak
functions used in his solution of the Christoffel problem, since they are closely related to
the action of the Hard Lefschetz integration operator on Minkowski valuations (see Section 5).
In the final part of this section we give a new proof of the bijectivity of integral transforms involving C.\ Berg's functions.
For the background material we refer the reader to \textbf{\cite[\textnormal{Chapter 8.3}]{schneider93}}, \textbf{\cite{Groemer1996}}, and \textbf{\cite{morimoto98}}.

We write $\Delta_S$ for the Laplacian (or Laplace--Beltrami
operator) on $S^{n-1}$. If $f, g \in C^2(S^{n-1})$, then we have
\[\int_{S^{n-1}} f(u)\,\Delta_Sg(u)\,du = \int_{S^{n-1}} g(u)\,\Delta_S f(u)\,du.  \]

The finite dimensional vector space of spherical harmonics of
dimension $n$ and degree $k$ will be denoted by $\mathcal{H}_k^n$
and we write $N(n,k)$ for its dimension. Spherical harmonics are eigenfunctions of $\Delta_S$, more precisely, for $Y_k \in
\mathcal{H}_k^n$,
\begin{equation} \label{deltasmult}
\Delta_S Y_k = -k(k + n - 2)\,Y_k.
\end{equation}

Let $L^2(S^{n-1})$ denote the Hilbert space of square-integrable
functions on $S^{n-1}$ with the usual inner product
$(\,\cdot\,,\,\cdot\,)$. The spaces $\mathcal{H}_k^n$ are
pairwise \linebreak orthogonal with respect to this inner product. If $\{Y_{k,1}, \ldots, Y_{k,N(n,k)}\}$ is an orthonormal basis of
$\mathcal{H}_k^n$, then the collection $\{Y_{k,1}, \ldots, Y_{k,N(n,k)}: k \in
\mathbb{N}\}$ is a complete orthogonal system in $L^2(S^{n-1})$,
that is, the Fourier series
\begin{equation} \label{fourierexp}
 f \sim \sum_{k=0}^{\infty} \pi_k f
\end{equation}
converges to $f$ in the $L^2$ norm for every $f \in L^2(S^{n-1})$. Here, we used $\pi_k: L^2(S^{n-1}) \rightarrow \mathcal{H}_k^n$ to denote the orthogonal projection.
Since the \emph{Legendre polynomial} $P_k^n \in C([-1,1])$ of dimension $n$ and degree $k$ satisfies
\begin{equation*}
\sum_{i=1}^{N(n,k)} Y_{k,i}(u)\,Y_{k,i}(v) = \frac{N(n,k)}{\omega_n}\, P_k^n(u \cdot v),
\end{equation*}
where $\omega_m$ denotes the surface area of the $m$-dimensional unit ball, we have
\begin{equation} \label{projleg}
(\pi_k f)(v) = \sum_{i=1}^{N(n,k)} (f,Y_{k,i})Y_{k,i}(v) = \frac{N(n,k)}{\omega_n} \int_{S^{n-1}} f(u)\, P_k^n(u\cdot v) \,du.
\end{equation}

Throughout the article we use $\bar{e} \in S^{n-1}$ to denote the pole of the sphere and we write $\mathrm{SO}(n-1)$ for the stabilizer
in $\mathrm{SO}(n)$ of $\bar{e}$. A function or measure on $S^{n-1}$ is called \emph{zonal} if it is $\mathrm{SO}(n-1)$ invariant. Clearly, zonal functions
depend only on the value of $u \cdot \bar{e}$.

The subspace of zonal functions in $\mathcal{H}_k^n$ is $1$-dimensional for every $k \in \mathbb{N}$ and spanned by the function
$u \mapsto P_k^n(u \cdot \bar{e})$. Since the spaces $\mathcal{H}_k^n$ are invariant under the natural action of $\mathrm{SO}(n)$,
the functions $u \mapsto P_k^n(u \cdot v)$, for fixed $v \in S^{n-1}$, are elements of $\mathcal{H}_k^n$. The orthogonality
of the spaces $\mathcal{H}_k^n$ is reflected by the fact that the Legendre polynomials $P_k^n$ form a complete orthogonal system with respect to the inner product on $C([-1,1])$ defined by
\[\left [ p,q \right ]_n = \int_{-1}^1 p(t)\,q(t)\,(1-t^2)^{\frac{n-3}{2}}\,dt.  \]
From the orthogonality property of the Legendre polynomials and (\ref{projleg}), it is not difficult to show that any function $\phi \in L^2([-1,1])$
(or, equivalently, any zonal $g \in L^2(S^{n-1})$) admits a series expansion
\begin{equation} \label{expzonal}
\phi \sim \sum_{k=0}^{\infty} \frac{N(n,k)}{\omega_n}\, a_k^n[\phi]\,P_k^n,
\end{equation}
where
\begin{equation} \label{multleg}
a_k^n[\phi] = \omega_{n-1} \int_{-1}^1 \phi(t)\,P_k^n(t)\,(1-t^2)^{\frac{n-3}{2}}\,dt= \omega_{n-1}\left [ P_k^n,\phi \right ]_n.
\end{equation}

For the explicit calculation of integrals of the form (\ref{multleg}) the following formula of Rodrigues for the Legendre polynomials
is often very useful:
\begin{equation} \label{formrod}
P_k^n(t) = \frac{(-1)^k}{2^k\left (\frac{n-1}{2}\right)_k} (1-t^2)^{-\frac{n-3}{2}} \frac{d^k}{dt^k} (1-t^2)^{\frac{n-3}{2}+k},
\end{equation}
where, for $\alpha \in \mathbb{R}$ and $k \in \mathbb{N}$, we have used $(\alpha)_k$ to abbreviate the product $\alpha(\alpha + 1) \cdots (\alpha + k - 1)$.
Using (\ref{formrod}) one can show that the derivatives of Legendre polynomials are again Legendre polynomials. For $l \geq k$, we have
\begin{equation} \label{ablleg}
\frac{d^k}{dt^k}P_l^n(t) = 2^k \left (\frac{n}{2} \right )_k \frac{N(n + 2k,l-k)}{N(n,l)}\, P_{l-k}^{n+2k}.
\end{equation}

Next we recall the \emph{Gegenbauer polynomials} which can be defined for $\alpha > 0$ by means of
the generating function
\[\frac{1}{(1+r^2-2rt)^{\alpha}} = \sum_{k=0}^{\infty} C_k^{\alpha}(t)\,r^n.  \]

\pagebreak

\noindent For $n \geq 3$, their relation to Legendre polynomials can be expressed by
\begin{equation} \label{gegenbleg}
C_k^{(n-2)/2} = {n + k - 3 \choose n - 3} P_k^n.
\end{equation}

For the following well known auxiliary result about the
spherical harmonic expansion of smooth functions, see, e.g.,
\textbf{\cite[\textnormal{p.\ 36}]{morimoto98}}.

\begin{lem} \label{smoothharmonic} If $f \in C^{\infty}(S^{n-1})$, then the sequence
$\|\pi_k f\|_{\infty}$, $k \in \mathbb{N}$, is \emph{rapidly
decreasing}, that is, for any $m \in \mathbb{N}$, we have $\sup
\{k^m\|\pi_k f\|_{\infty}: k \in \mathbb{N}\} < \infty.$
Conversely, if $Y_k \in \mathcal{H}_k^n$, $k \in \mathbb{N}$, is
a sequence of spherical harmonics such that $\|Y_k\|_{\infty}$ is
rapidly decreasing, then the function
\[f(u) = \sum_{k=0}^\infty Y_k(u), \qquad u \in S^{n-1},  \]
is $C^{\infty}$ and $\pi_k f = Y_k$ for every $k \in \mathbb{N}$.
\end{lem}

For $f \in C^{\infty}(S^{n-1})$ and $m \in \mathbb{N}$, define
\[(-\Delta_S)^{\frac{m}{2}}f = \sum_{k=0}^{\infty} (k(k+n-2))^{\frac{m}{2}} \pi_k f.  \]
Note that, by Lemma \ref{smoothharmonic},
$(-\Delta_S)^{\frac{m}{2}}f \in C^{\infty}(S^{n-1})$.

If we endow the vector space $C^{\infty}(S^{n-1})$ with the
topology defined by the family of semi norms
$\|(-\Delta_S)^{\frac{m}{2}}f\|_{\infty}$, $m \in \mathbb{N}$,
then $C^{\infty}(S^{n-1})$ becomes a Fr\'echet space. Moreover,
the spherical harmonic expansion (\ref{fourierexp}) of any $f \in
C^{\infty}(S^{n-1})$ converges to $f$ in this topology.

A \emph{distribution} on $S^{n-1}$ is a continuous linear
functional on $C^{\infty}(S^{n-1})$. We write
$C^{-\infty}(S^{n-1})$ for the space of distributions on
$S^{n-1}$ equipped with the topology of weak convergence and use
$\langle \,\cdot\,,\,\cdot\, \rangle$ to denote the canonical
bilinear pairing on $C^{\infty}(S^{n-1}) \times
C^{-\infty}(S^{n-1})$.

A (signed) measure $\sigma$ on $S^{n-1}$ defines a distribution
$T_{\sigma}$ by
\[\langle f, T_{\sigma} \rangle = \int_{S^{n-1}} f(u)\,d\sigma(u), \qquad f \in C^{\infty}(S^{n-1}).  \]
Using the continuous linear injection $\sigma \mapsto T_{\sigma}$, we
can regard $\mathcal{M}(S^{n-1})$ as a subspace of
$C^{-\infty}(S^{n-1})$. In the same way, the spaces
$C^{\infty}(S^{n-1})$, $C(S^{n-1})$, and $L^2(S^{n-1})$ can be
viewed as subspaces of $C^{-\infty}(S^{n-1})$ and we have
\begin{equation} \label{inclusion1742}
C^{\infty}(S^{n-1}) \subseteq C(S^{n-1}) \subseteq L^2(S^{n-1})
\subseteq \mathcal{M}(S^{n-1}) \subseteq C^{-\infty}(S^{n-1}).
\end{equation}

Since $\pi_k: L^2(S^{n-1}) \rightarrow
\mathcal{H}_k^n$ is self-adjoint,
that is, $(\pi_k f,g) = (f,\pi_k g)$ for all $f, g \in
L^2(S^{n-1})$ and $k \in \mathbb{N}$, we define the $k$-spherical harmonic component
$\pi_k T$ of a distribution $T \in C^{-\infty}(S^{n-1})$ as the
distribution given by
\[\langle f, \pi_k T \rangle = \langle \pi_k f, T \rangle, \qquad f \in C^{\infty}(S^{n-1}).  \]

\begin{lem} \label{distharmonic} {\bf ($\!\!$\cite[\textnormal{p.\ 38}]{morimoto98})} If $T \in C^{-\infty}(S^{n-1})$, then $\pi_k T \in
\mathcal{H}_k^n$ for every $k \in \mathbb{N}$ and the sequence
$\|\pi_k T\|_{\infty}$, $k \in \mathbb{N}$, is \emph{slowly
increasing}, that is, there exist $C > 0$ and $j \in \mathbb{N}$
such that $\|\pi_k T\|_{\infty} \leq C(1 + k^j)$ for every $k \in
\mathbb{N}$.\\
Conversely, if $Y_k \in \mathcal{H}_k^n$, $k \in \mathbb{N}$, is
a sequence of spherical harmonics such that $\|Y_k\|_{\infty}$ is
slowly increasing, then
\[\langle g, T \rangle  = \sum_{k=0}^\infty \int_{S^{n-1}} g(u)Y_k(u) \,du, \qquad g \in C^{\infty}(S^{n-1}),  \]
defines a distribution $T \in C^{-\infty}(S^{n-1})$ for which
$\pi_k T = Y_k$ for every $k \in \mathbb{N}$.
\end{lem}

We can also extend the Laplacian
to distributions $T \in C^{-\infty}(S^{n-1})$, by defining $\Delta_S T$ as the distribution given by
\[\langle f,\Delta_S T \rangle = \langle \Delta_S f, T \rangle, \qquad f \in C^{\infty}(S^{n-1}).  \]
Note that, by (\ref{inclusion1742}), $\Delta_S$ can now also act
on continuous functions on $S^{n-1}$. This is of particular
importance for us, since the support function $h(K,\cdot)$ and
the first-order area measure $S_1(K,\cdot)$ of a convex body $K
\in \mathcal{K}^n$ are related \nolinebreak by
\begin{equation} \label{boxhks1}
\Box_n h(K,\cdot) = S_1(K,\cdot),
\end{equation}
where $\Box_n$ is the differential operator given by
\[\Box_n h = h +  \frac{1}{n-1}\Delta_Sh. \]

From the definition of $\Box_n$ and (\ref{deltasmult}), we see
that for $f \in C^{\infty}(S^{n-1})$ the spherical harmonic
expansion of $\Box_n f$ is given by
\begin{equation} \label{boxnmult}
\Box_n f \sim \sum_{k=0}^{\infty} \frac{(1-k)(k+n-1)}{n-1} \pi_kf.
\end{equation}
Thus, the kernel of the linear operator $\Box_n:
C^{\infty}(S^{n-1}) \rightarrow C^{\infty}(S^{n-1})$ is given by
$\mathcal{H}_1^n$ and consists precisely of the
restrictions of linear functions on $\mathbb{R}^n$ to $S^{n-1}$.
Let $C_{\mathrm{o}}^{\infty}(S^{n-1})$ denote the Fr\'echet
subspace of $C^{\infty}(S^{n-1})$ given by
\[C_{\mathrm{o}}^{\infty}(S^{n-1}) = \{f \in C^{\infty}(S^{n-1}): \pi_1 f = 0\}  \]
and define $C_{\mathrm{o}}^{-\infty}(S^{n-1})$ analogously.

\pagebreak

Since the linear operator $\Box_n:
C_{\mathrm{o}}^{\infty}(S^{n-1}) \rightarrow
C_{\mathrm{o}}^{\infty}(S^{n-1})$ is an isomorphism, it is a
natural problem to find an (explicit) inversion formula. This was
accomplished by C.\ Berg \textbf{\cite{cberg}} in the late 1960s and, due to
(\ref{boxhks1}), is closely related to his solution of the
classical Christoffel problem which consists in finding necessary
and sufficient conditions for a Borel measure on $S^{n-1}$ to be
the first-order area measure of a convex body.

In order to describe C.\ Berg's inversion formula for $\Box_n$, let us recall
the \emph{Funk--Hecke Theorem}: If $\phi \in C([-1,1])$ and
$\mathrm{F}_{\phi}$ is the integral transform on
$\mathcal{M}(S^{n-1})$ defined by
\[(\mathrm{F}_{\phi}\sigma)(u) = \int_{S^{n-1}} \phi(u \cdot v)\,d\sigma(v), \qquad u \in S^{n-1},  \]
then the spherical harmonic expansion
of $\mathrm{F}_{\phi}\sigma \in C(S^{n-1})$ is given by
\begin{equation} \label{funkhecke}
\mathrm{F}_{\phi}\sigma \sim \sum_{k=0}^\infty
a_k^n[\phi]\,\pi_k\sigma,
\end{equation}
where the numbers $a_k^n[\phi]$ are given by (\ref{multleg}) and called the \emph{multipliers}
of $\mathrm{F}_{\phi}$.

Using the theory of subharmonic functions on $S^{n-1}$, C.\ Berg
proved that for every $n \geq 2$ there exists a uniquely
determined $C^\infty$ function $g_n$ on $(-1,1)$ such that the zonal function $u \mapsto g_n(u \cdot \bar{e})$ is in $L^1(S^{n-1})$ and
\begin{equation} \label{multgn}
a_1^n[g_n] = 0, \quad \qquad a_k^n[g_n] = \frac{n-1}{(1-k)(k+n-1)}, \quad k \neq 1.
\end{equation}
For later reference, we just state here
\begin{equation} \label{g2}
g_2(t) = \frac{1}{2\pi} \left( (\pi-\arccos t )(1-t^2)^{\frac{1}{2}} - \frac{t}{2} \right)
\end{equation}
and
\begin{equation} \label{g3}
g_3(t) = \frac{1}{2\pi} \left( 1 + t \ln(1-t) + \left( \frac 4 3 - \ln 2 \right) t \right ).
\end{equation}
We note that, by (\ref{multgn}), our normalization of the $g_n$ differs from C.\ Berg's original one.
It follows from (\ref{boxnmult}), (\ref{funkhecke}), and (\ref{multgn}) that
\[f(u) = \int_{S^{n-1}} g_n(u \cdot v)(\Box_nf)(v)\,dv, \qquad u \in S^{n-1},  \]
for every $f \in C_{\mathrm{o}}^{\infty}(S^{n-1})$, which is the
desired inversion formula. However, for our purposes we need the following more general fact.

\begin{theorem} \label{fgjiso} For every $n \geq 2$ and $2 \leq j \leq n$, the integral transform $\mathrm{F}_{g_j}:
C_{\mathrm{o}}^{\infty}(S^{n-1}) \rightarrow
C_{\mathrm{o}}^{\infty}(S^{n-1})$, given by
\[(\mathrm{F}_{g_j}f)(u) = \int_{S^{n-1}} g_j(u \cdot v) f(v)\,dv, \qquad u \in S^{n-1},  \]
is an isomorphism.
\end{theorem}

Theorem \ref{fgjiso} follows for example from a recent result of Goodey and Weil \textbf{\cite[\textnormal{Theorem 4.3}]{GoodeyWeil2014}}. However, we give a different and more elementary proof
below that also yields additional information of independent interest. For this, note that, by Lemma \ref{smoothharmonic}, it is sufficient
to show that the multipliers $a_k^n[g_j]$ are non-zero for $k \neq 1$ and that they are slowly increasing. Therefore,
Theorem \ref{fgjiso} is a direct consequence of the following.

\begin{theorem} \label{multipliers} For $n \geq 2$, $2 \leq j \leq n$, and $k \neq 1$, we have
\[a_k^n[g_j] = -\frac{\pi^{\frac{n-j}{2}}(j-1)}{4}\frac{\Gamma\left(\frac{n-j+2}{2}\right)\Gamma \left ( \frac{k-1}{2} \right )
 \Gamma \left (\frac{j+k-1}{2} \right )}{\Gamma \left (\frac{n-j+k+1}{2}\right ) \Gamma \left (\frac{n+k+1}{2}\right )}.\]
\end{theorem}
{\it Proof.}  For $n \geq 2$, $d \geq 0$, and $k \neq 1$, by (\ref{multleg}), we have to determine
\begin{equation}\label{eq: multipliers definition}
a^{n,d}_k := a_k^{n+d}[g_n] = \omega_{n+d-1} \left [ P_k^{n+d},g_n \right ]_{n+d},
\end{equation}
where we know from (\ref{multgn}) that
\begin{equation}\label{eq: multipliers d=0}
a^{n,0}_k = \frac {n-1} {(1-k)(n-1+k)} .
\end{equation}

We start with the case $d = 1$. By (\ref{expzonal}) and (\ref{multgn}), we have
\begin{equation*}
g_n = \sum_{l=0}^\infty \frac{N(n,l)}{\omega_n}\, a^{n,0}_l P^n_l,
\end{equation*}
where the sum converges in the topology induced by $[\,\cdot\,,\cdot\,]_n$, which implies the convergence in
the topology induced by $[\,\cdot\,,\cdot\,]_{n+1}$. Consequently,
\begin{equation} \label{an1k}
a^{n,1}_k = \omega_n \left [ P^{n+1}_k , g_n \right ]_{n+1} = \sum_{l=0}^\infty N(n,l)\, a^{n,0}_l \left [ P^{n+1}_k, P^n_l \right ]_{n+1}.
\end{equation}
Since Legendre polynomials of degree $k$ are even if $k$ is even, and odd otherwise, we may assume that $k$ and $l$ have the same parity.
Since $\left [P^{n+1}_k, P^n_l \right ]_{n+1}$ vanishes for $l < k$ (see the next calculation), let $l \geq k$ and put
\[\beta :=  \frac{n-2}{2}. \]

\pagebreak

\noindent If $\beta + k \geq \frac{1}{2}$, that is, $(n,k) \neq (2,0)$, then it follows from (\ref{formrod}), integration by parts, (\ref{ablleg}), and (\ref{gegenbleg}) that
\begin{align*}
\left [ P^{n+1}_k, P^n_l \right ]_{n+1}
        & = \frac{(-1)^k}{2^k (\beta + 1)_k} \int_{-1}^1 \left( \frac{d^k}{dt^k} (1-t^2)^{\beta + k} \right)P^n_l(t)\,dt \\
        & = \frac{1}{2^k (\beta + 1)_k} \int_{-1}^1 (1-t^2)^{\beta + k} \left( \frac{d^k}{dt^k} P^n_l(t) \right) dt \\
        & = \frac{N(n+2k, l-k)}{N(n,l)} \int_{-1}^1 (1-t^2)^{\beta + k}\, P^{n+2k}_{l-k}(t)\,dt\\
        & = \frac{\beta + l}{N(n,l)(\beta + k)} \int_{-1}^1 (1-t^2)^{\beta + k}\, C^{\beta + k}_{l-k}(t)\,dt .
    \end{align*}
For $\alpha \in \frac{1}{2} \mathbb{N}$ and even $m$, we have (cf.\ \textbf{\cite[\textnormal{p. 424}]{GoodeyWeil1992}})
\begin{equation*}
c^\alpha_m = \int_{-1}^1 (1-t^2)^\alpha\, C^\alpha_m(t)\,dt =
        -\frac {\alpha\, 4^{\alpha + \frac{1}{2}}\, m!\, \Gamma\left (\frac{m}{2} + \alpha + 1 \right )^2}
        {(m-1) \left (\frac{m}{2} + \alpha \right ) (m + 2 \alpha + 1)!\, \Gamma\left (\frac{m}{2} + 1\right )^2}.
\end{equation*}
Plugging this into (\ref{an1k}) and changing the summation index, yields
\begin{equation*}
a^{n,1}_k = \sum_{l=0}^\infty \frac {\beta+k+2l} {\beta+k}\, a^{n,0}_{k+2l}\, c^{\beta+k}_{2l} =
        \left( \beta + \frac 1 2 \right) 4^{\beta+k+1} \sum_{l=0}^\infty q(\beta, k, l) ,
\end{equation*}
where
\begin{equation*}
q(\beta, k, l) = \frac{2l\,(2l{-}2)!\,(\beta{+}k{+}2l)\,\Gamma(\beta{+}k{+}l{+}1)^2}
{(k{+}2l{-}1)\, (2\beta{+}k{+}2l{+}1)\, (\beta{+}k{+}l)\, (2\beta{+}2k{+}2l{+}1)!\, \Gamma(l{+}1)^2} .
\end{equation*}

Using Zeilberger's algorithm (see, e.g., \textbf{\cite{PetkovsekWilfZeilberger1996}}), we find that $q$ satisfies the following recurrence relation
\begin{equation*}
A(\beta,\! k) q(\beta+1,\! k,\! l) + B(\beta,\! k) q(\beta,\! k,\! l) = q(\beta,\! k,\! l+1) C(\beta,\! k,\! l+1) - q(\beta,\! k,\! l) C(\beta,\! k,\! l) ,
\end{equation*}
where
\begin{equation*}
A(\beta,\! k) = 4 (2\beta{+}k{+}4), \quad B(\beta,\! k) = {-}(2\beta{+}k{+}1), \quad C(\beta,\! k,\! l) = {-}\frac {l (k{+}2l{-}1)} {\beta{+}k{+}2l} .
\end{equation*}
If we let $Q(\beta, k) = \sum_{l=0}^\infty q(\beta, k,l)$, then we obtain
\begin{equation*}
Q(\beta+1, k) = \frac {2\beta+k+1} {4(2\beta+k+4)}\, Q(\beta, k)
\end{equation*}
or, in terms of the multipliers,
\begin{equation}\label{eq: recurrence d=1, one}
a^{n+2,1}_k = \frac {(n+1)\,(n+k-1)}{(n-1)\,(n+k+2)}\, a^{n,1}_k.
\end{equation}

\pagebreak

The function $q$ also satisfies the recurrence relation
\begin{equation*}
D(\beta,\! k) q(\beta,\! k+2,\! l) + E(\beta,\! k) q(\beta,\! k,\! l) = q(\beta,\! k,\! l+1) F(\beta,\! k,\! l+1) - q(\beta,\! k,\! l) F(\beta,\! k,\! l),
\end{equation*}
where
\begin{equation*}
D(\beta, k) = 16 (k+2) (2\beta+k+4) ,\quad E(\beta, k) = -(k-1) (2\beta+k+1)
\end{equation*}
and
\[F(\beta, k ,l) = - \frac {1}{(\beta+k+2l) (2\beta+2k+2l+3)} \sum_{i=1}^3 l^i\, p_i(\beta,k),\]
with polynomials $p_1, p_2, p_3$ given by
\begin{align*}
p_1(\beta,k) & = 8\beta^3 + 16k\beta^2 + 20\beta^2 + 12k^2\beta + 26k\beta + 10\beta + 4k^3 + 9k^2 + 4k + 3, \\
p_2(\beta,k) & = 16 \beta^2 + 24k\beta + 32\beta + 12k^2 + 24k + 4, \\
p_3(\beta,k) & = 8\beta + 8k + 12.
\end{align*}
Summing again over all $l$, we arrive at
\begin{equation*}
Q(\beta,k+2) = \frac{(k-1)\,(2\beta+k+1)} {16(k+2)\,(2\beta+k+4)}\, Q(\beta, k) .
\end{equation*}
In terms of the multipliers this means
\begin{equation}\label{eq: recurrence d=1, two}
a^{n,1}_{k+2} = \frac {(k-1)\,(n+k-1)} {(k+2)\,(n+k+2)}\, a^{n,1}_k .
\end{equation}

In order to solve \eqref{eq: recurrence d=1, one} and \eqref{eq: recurrence d=1, two}, we need four initial values of $a^{n,1}_k$.
We also have to calculate $a^{2,1}_0$, which was not covered by the above arguments. Using (\ref{g2}), (\ref{g3}), and \eqref{eq: multipliers definition}, elementary integration
yields
\begin{equation}\label{eq: initial values d=1, one}
a^{2,1}_0 = \frac{\pi^2}{4}, \quad a^{2,1}_2 = -\frac{\pi^2}{32}, \quad a^{2,1}_3 = -\frac{4}{45}, \quad
a^{3,1}_0 = \frac {2\pi}{3}, \quad a^{3,1}_3 = -\frac{\pi}{24}.
\end{equation}
This leads to the sequence
\begin{equation}\label{eq: multipliers d=1}
a^{n,1}_k = - \frac{\pi}{8} (n-1) \frac{\Gamma\left( \frac{k-1}{2} \right) \Gamma\left( \frac {n+k-1}{2} \right)}
{\Gamma\left( \frac {k+2} 2 \right) \Gamma\left( \frac {n+k+2} 2 \right)}
\end{equation}
which satisfies \eqref{eq: recurrence d=1, one}, \eqref{eq: recurrence d=1, two},
and has the initial values \eqref{eq: initial values d=1, one}.

Now let $d \geq 0$ be arbitrary. For $l \geq 2$, the Legendre polynomials satisfy the recurrence relation
(see, e.g., \textbf{\cite[\textnormal{Lemma 3.3.10}]{Groemer1996}})
\[(n+d+2l-2)(n+d-1) P^{n+d}_l = (n+d+l-2)(n+d+l-1) P^{n+d+2}_l - (l-1)l P^{n+d+2}_{l-2}.\]

\pagebreak

\noindent From this and the fact that $P_0^m(t) = 1$, $P^m_1(t) = t$ for all $m \geq 1$,
we obtain
\begin{align*}
a^{n,d+2}_k & = \omega_{n+d+1} \left [ P^{n+d+2}_k,g_n \right ]_{n+d+2} \\
            & = \omega_{n+d+1} \sum_{l=0}^\infty \frac{N(n+d,l)}{\omega_{n+d}} \, a^{n,d}_l \left [ P^{n+d+2}_k,P^{n+d}_l \right ]_{n+d+2} \\
            & = \frac{2\pi}{n+d+2k} \left( a^{n,d}_k - a^{n,d}_{k+2} \right ).
\end{align*}
Finally, the sequence which solves this recurrence relation and has the initial values \eqref{eq: multipliers d=0} and \eqref{eq: multipliers d=1}
is given by
\begin{equation*}
a^{n,d}_k = - \frac {\pi^{\frac{d}{2}}\,(n-1)}{4} \, \frac{\Gamma\left( \frac {d+2} 2 \right) \Gamma\left( \frac {k-1} 2 \right) \Gamma\left( \frac {n+k-1} 2 \right)}
{\Gamma\left( \frac {d+k+1} 2 \right) \Gamma\left( \frac {n+d+k+1} 2 \right)} ,\qquad k \neq 1.
\end{equation*}

\vspace{-0.4cm}

\hfill $\blacksquare$

\vspace{0.6cm}

We end this section with the following important definition, given rise to by Theorem \ref{fgjiso}.

\vspace{0.3cm}

\noindent {\bf Definition} \emph{For $2 \leq j \leq n$, let $\Box_j: C_{\mathrm{o}}^{\infty}(S^{n-1}) \rightarrow
C_{\mathrm{o}}^{\infty}(S^{n-1}) $ denote the linear operator which is
inverse to the integral transform $\mathrm{F}_{g_j}$.}

\vspace{1cm}

\centerline{\large{\bf{ \setcounter{abschnitt}{5}
\arabic{abschnitt}. Generalized valuations and Minkowski
valuations}}}

\reseteqn \alpheqn \setcounter{theorem}{0}

\vspace{0.6cm}

In the following we recall several results on translation invariant
(scalar and convex body valued) valuations, in
particular, the product structure on smooth valuations and the Alesker--Poincar\'e duality.
We also discuss basic properties of the Hard Lefschetz operators and a new isomorphism between
generalized valuations of degree one and generalized functions on the sphere. At the end of this section,
we state a recent representation theorem for Minkowski valuations intertwining rigid motions and give an alternative
description of the classes $\mathbf{MVal}_{i,j}^{\mathrm{SO}(n)}$.

A map $\mu$ defined on convex bodies in $\mathbb{R}^n$ and taking
values in an Abelian semigroup $A$ is called a {\it valuation} or
{\it additive} if
\[\mu(K) + \mu(L) = \mu(K \cup L) + \mu(K \cap L)   \]
whenever $K \cup L$ is convex. If $G$ is a group of affine
transformations on $\mathbb{R}^n$, a valuation $\mu$ is called
$G$-invariant if $\mu(gK) = \mu(K)$ for all $K \in \mathcal{K}^n$
and $g \in G$.

Let $\mathbf{Val}$ denote the vector space of continuous translation invariant
scalar valued valuations. The structure theory of translation invariant
valuations has its starting point in a classical result of McMullen
\textbf{\cite{McMullen77}}, who showed that
\begin{equation} \label{mcmullen}
\mathbf{Val} = \bigoplus_{0\leq i \leq n} \mathbf{Val}_i^+ \oplus
\mathbf{Val}_i^-,
\end{equation}
where $\mathbf{Val}_i^+ \subseteq \mathbf{Val}$ denotes the
subspace of {\it even} valuations (homogeneous) of degree $i$, and
$\mathbf{Val}_i^-$ denotes the subspace of {\it odd} valuations of
degree $i$. The
space $\mathbf{Val}$ becomes a Banach space, when endowed with
the norm
\[\|\mu \|=\sup \{|\mu(K)|: K \subseteq B\}.  \]
The general linear group $\mathrm{GL}(n)$ acts on the Banach
space $\mathbf{Val}$ in a natural way: For every $A \in
\mathrm{GL}(n)$ and $K \in \mathcal{K}^n$,
\[(A\cdot\mu)(K)=\mu(A^{-1}K), \qquad \mu \in \mathbf{Val}.  \]
Note that the subspaces $\mathbf{Val}_i^\pm$ are invariant under this $\mathrm{GL}(n)$-action.
In fact, a deep result of Alesker \textbf{\cite{Alesker01}}, known as the
Irreducibility Theorem, states that these subspaces are also irreducible:

\begin{theorem} \label{irr} \emph{(Alesker \textbf{\cite{Alesker01}})} The
natural representation of $\mathrm{GL}(n)$ on
$\mathbf{Val}_i^{\pm}$ is irreducible for any $i \in \{0, \ldots,
n\}$.
\end{theorem}

It follows from Theorem \ref{irr} that any
$\mathrm{GL}(n)$-invariant subspace of translation invariant
continuous valuations (of a given degree $i$ and parity) is already
dense in $\mathbf{Val}_i^\pm$.

\vspace{0.3cm}

\noindent {\bf Definition} \emph{A valuation $\mu \in
\mathbf{Val}$ is called smooth if the map $\mathrm{GL}(n)
\rightarrow \mathbf{Val}$ defined by $A \mapsto A\cdot\mu$ is
infinitely differentiable.}

\vspace{0.3cm}

The subspace of smooth translation invariant valuations is
denoted by $\mathbf{Val}^{\infty}$, and we write
$\mathbf{Val}_i^{\pm,\infty}$ for smooth valuations in
$\mathbf{Val}_i^{\pm}$. It is well known (cf.\
\textbf{\cite[\textnormal{p.\ 32}]{wallach1}}) that $\mathbf{Val}_i^{\pm,\infty}$ is a dense
$\mathrm{GL}(n)$ invariant subspace of $\mathbf{Val}_i^{\pm}$.
Moreover, $\mathbf{Val}^{\infty}$ carries a natural Fr\'echet
space topology, called G\aa rding topology (see
\textbf{\cite[\textnormal{p.\ 33}]{wallach1}}), which is stronger
than the topology induced from $\mathbf{Val}$. Finally, we note
that the representation of $\mathrm{GL}(n)$ on
$\mathbf{Val}^{\infty}$ is continuous.

\vspace{0.3cm}

\noindent {\bf Examples:}
\begin{enumerate}
\item[(a)] If $L \in \mathcal{K}^n$ is strictly convex with
smooth boundary, then
\[ \mu_L: \mathcal{K}^n \rightarrow \mathbb{R}, \qquad \mu_L(K)=V_n(K+L),  \]
is a smooth valuation.

\pagebreak

\item[(b)] If $f \in C^{\infty}_{\mathrm{o}}(S^{n-1})$ and $0 \leq i \leq n - 1$, then $\nu_{i,f}: \mathcal{K}^n \rightarrow
\mathbb{R}$ defined by
\begin{equation} \label{lesnope1}
\nu_{i,f}(K)= \int_{S^{n-1}} f(u)\,dS_i(K,u),
\end{equation}
is a smooth valuation in $\mathbf{Val}_i^{\infty}$.
\end{enumerate}

\vspace{0.2cm}

Before we turn to generalized valuations, we recall the definition of the Alesker product of
smooth translation invariant valuations.

\begin{theorem} {\bf (\!\!\cite{Alesker04a})} There exists a
bilinear product
\[\mathbf{Val}^{\infty} \times \mathbf{Val}^{\infty} \rightarrow
\mathbf{Val}^{\infty}, \quad (\mu,\nu) \mapsto \mu \cdot
\nu,\] which is uniquely determined by the following two
properties:
\begin{enumerate}
\item[(i)] The product is continuous in the G\aa rding topology.
\item[(ii)] If $L_1, L_2 \in \mathcal{K}^n$ are strictly convex and
smooth, then
\[(\mu_{L_1} \cdot \mu_{L_2})(K) = V_{2n}(\iota(K)+ L_1 \times L_2),  \]
where $\iota: \mathbb{R}^n \rightarrow \mathbb{R}^n \times
\mathbb{R}^n$ is defined by $\iota(x)=(x,x)$.
\end{enumerate}
Endowed with this multiplicative structure,
$\mathbf{Val}^{\infty}$ becomes an associative and commutative
algebra which is graded by the degree of homogeneity and with unit
given by the Euler characteristic.
\end{theorem}

\vspace{0.1cm}

The next example was computed in \textbf{\cite{Alesker04a}} and will be needed in the appendix.

\vspace{0.3cm}

\noindent {\bf Example:}

\vspace{0.1cm}

\noindent Let $L_1, \ldots, L_{n-i} \in \mathcal{K}^n$ and $M_1, \ldots, M_i \in \mathcal{K}^n$ be strictly convex and smooth. If
$\mu \in \mathbf{Val}_i^{\infty}$ and $\nu \in \mathbf{Val}_{n-i}^{\infty}$ are defined by
\[\mu(K)=V(K[i],L_1,\ldots,L_{n-i}) \quad \mbox{and} \quad \nu(K)=V(K[n-i],M_1,\ldots,M_i),  \]
then
\begin{equation} \label{specialprod}
(\mu \cdot \nu)(K) = {n \choose
i}^{-1}V(-L_1,\ldots,-L_{n-i},M_1,\ldots,M_i)\,V_n(K).
\end{equation}

\vspace{0.2cm}

The above example is just a special case of the more general fact that the Alesker product gives
rise to a nondegenerate bilinear pairing between smooth valuations of complementary degree.

\begin{theorem} \label{poincaredual} {\bf (\!\!\cite{Alesker04a})} For every $0 \leq i
\leq n$, the continuous bilinear pairing
\[<\!\cdot\,,\cdot\!>\,:\mathbf{Val}^{\infty}_i \times \mathbf{Val}^{\infty}_{n-i} \rightarrow
\mathbf{Val}_n, \qquad (\mu,\nu) \mapsto \mu \cdot \nu,\] is
nondegenerate. In particular, the induced Poincar\'e duality map
\[\mathbf{Val}^{\infty}_i \rightarrow \left (\mathbf{Val}^{\infty}_{n-i}\right )^* \otimes \mathbf{Val}_n, \qquad \mu \mapsto\, <\!\mu,\cdot>, \]
is continuous, injective and has dense image with respect to the
weak topology.
\end{theorem}

Here and in the following, for a Fr\'echet space $X$, we denote by $X^*$ its topological dual endowed with the weak topology.

Motivated by Theorem \ref{poincaredual}, the notion of generalized valuations was introduced recently \textbf{\cite{AleskerFaifman}}.
Before we state the definition, recall that by a classical theorem of Hadwiger \textbf{\cite[\textnormal{p.\ 79}]{hadwiger51}}
the space $\mathbf{Val}_n$ is spanned by the ordinary volume $V_n$. In other words, if we do not refer to any
Euclidean structure, then $\mathbf{Val}_n \cong \mathscr{D}(V)$, where $\mathscr{D}(V)$ denotes the vector space of all densities
on an $n$-dimensional vector space $V$ (see the Appendix for details).

\vspace{0.3cm}

\noindent {\bf Definition} \emph{The space of \emph{generalized valuations} is defined by
\[ \mathbf{Val}^{-\infty} = \left (\mathbf{Val}^{\infty} \right )^* \otimes \mathscr{D}(V) \]
and we define the space of generalized valuations of degree $i \in \{0,\ldots, n\}$ by}
\[ \mathbf{Val}_i^{-\infty} = \left (\mathbf{Val}_{n-i}^{\infty} \right )^* \otimes \mathscr{D}(V).  \]

\vspace{0.3cm}

By Theorem \ref{poincaredual}, we have a canonical embedding with dense image
\[\mathbf{Val}^{\infty} \hookrightarrow \mathbf{Val}^{-\infty}.  \]
Thus, $\mathbf{Val}^{-\infty}$ can be seen as a completion of $\mathbf{Val}^{\infty}$ in the weak topology.

\vspace{0.2cm}

In order to establish Theorem \ref{thm3}, we need the following new classification of generalized valuations of degree $1$.
A proof of this theorem was given by Semyon Alesker and is included in the appendix.

\begin{theorem} \label{genvalchar} The map
\[C^{\infty}_{\mathrm{o}}(S^{n-1}) \rightarrow
\mathbf{Val}_1^{\infty}, \qquad f \mapsto \left (K \mapsto \int_{S^{n-1}}f(u)\,h(K,u)\,du \right),  \] is an isomorphism which extends
uniquely by continuity in the weak topologies to an isomorphism
\[C^{-\infty}_{\mathrm{o}}(S^{n-1}) \rightarrow \mathbf{Val}_1^{-\infty}.  \]
\end{theorem}

Note that, by Theorem \ref{genvalchar}, if $\gamma \in \mathbf{Val}_1^{-\infty}$ and $T_{\gamma} \in C^{-\infty}_{\mathrm{o}}(S^{n-1})$
is the corresponding distribution, then we can evaluate $\gamma$ on convex bodies $K \in \mathcal{K}^n$ with smooth support function by
\begin{equation*}
\gamma(K) := \langle h(K,\cdot), T_{\gamma} \rangle.
\end{equation*}

\vspace{0.1cm}

Next we briefly recall the Hard Lefschetz operators on smooth translation invariant scalar valuations.
It is well known that McMullen's decomposition (\ref{mcmullen}) of the space $\mathbf{Val}$ implies a general Steiner type formula
for continuous translation invariant valuations which, in turn, gives rise to a derivation operator $\Lambda: \mathbf{Val} \rightarrow \mathbf{Val}$ defined by
\[(\Lambda \mu)(K) = \left . \frac{d}{dt} \right |_{t=0} \mu(K + tB).  \]
Note that $\Lambda$ commutes with the action of $\mathrm{O}(n)$ and that it preserves parity.
Moreover, if $\mu \in \mathbf{Val}_i$, then $\Lambda \mu \in \mathbf{Val}_{i-1}$.

The importance of the operator $\Lambda$ became evident from
a Hard Lefschetz type theorem established by Alesker \textbf{\cite{Alesker03}} for even
valuations and by Bernig and Br\"ocker \textbf{\cite{Bernig07b}} for general
valuations.
More recently, a dual version of this fundamental result was established in \textbf{\cite{Alesker04, alesker10}}. There, the
derivation operator $\Lambda$ is replaced by an integration operator $\mathfrak{L}: \mathbf{Val} \rightarrow \mathbf{Val}$
defined by
\begin{equation} \label{hardlefint}
(\mathfrak{L}\mu)(K) = (V_1\cdot \mu)(K) = \int_{\mathrm{AGr}_{n-1,n}} \mu(K \cap E)\,dE,
\end{equation}
where here and in the following $\mathrm{AGr}_{k,n}$ denotes the affine Grassmannian of $k$ planes in $\mathbb{R}^n$ and integration is
with respect to a (suitably normalized) invariant measure. The original definition of $\mathfrak{L}$
corresponds to the first equality in (\ref{hardlefint}) and it was proved by Bernig \textbf{\cite{Bernig07}} that the second equality holds. We also note that $\mathfrak{L}$ commutes with the action of $\mathrm{O}(n)$ and that it preserves parity.
Moreover, if $\mu \in \mathbf{Val}_i$, then $\mathfrak{L}\mu \in \mathbf{Val}_{i+1}$.

\vspace{0.2cm}

In the final part of this section we turn to Minkowski valuations. Recall that the trivial Minkowski valuation
maps every convex body to the set containing only the origin and that, for $0 \leq j \leq n$,
we denote by $\mathbf{MVal}_j$ the set of all continuous, translation invariant and $\mathrm{SO}(n)$ equivariant
Minkowski valuations of degree $j$. In the next lemma we state basic properties of such Minkowski
valuations which are well known (cf.\ \textbf{\cite{ABS2011, papschu12, Schu09}}) and are needed in what follows.

\pagebreak

\begin{lem} \label{basicMinkVal} If $\Phi_j \in \mathbf{MVal}_j$, $0 \leq j \leq n$, then the following
statements hold:
\begin{enumerate}
\item[(a)] The Steiner point of $\Phi_jK$ is at the origin, that is, $s(\Phi_jK) = o$,
for every $K \in \mathcal{K}^n$.
\item[(b)] There exists $r_{\Phi_j} \geq 0$ such that
\[W_{n-1}(\Phi_jK) = r_{\Phi_j}\, W_{n-j}(K)\]
for every $K \in \mathcal{K}^n$. If $\Phi_j$ is non-trivial, then $r_{\Phi_j} > 0$.
\item[(c)] The $\mathrm{SO}(n-1)$ invariant valuation $\nu_j \in \mathbf{Val}_j$, defined by
\[\nu_j(K) = h(\Phi_jK,\bar{e}) \]
uniquely determines $\Phi_j$ and is called \emph{the associated real valued valuation of}
$\Phi_j \in \mathbf{MVal}_j$.
\end{enumerate}
\end{lem}

Lemma \ref{basicMinkVal} (c) motivated the following definition which first appeared in \textbf{\cite{Schu09}}.

\vspace{0.3cm}

\noindent {\bf Definition} \emph{A Minkowski valuation $\Phi_j \in \mathbf{MVal}_j$, $0 \leq j \leq n$, is called \emph{smooth} if its associated real valued valuation $\nu_j$ is smooth.}

\vspace{0.3cm}

Recall that smooth translation invariant scalar valuations are dense in all continuous translation invariant scalar valuations. However, this does not
\emph{directly} imply the same for Minkowski valuations but instead additional arguments were needed for the proof which was given in \textbf{\cite{Schu09}} for even
and in \textbf{\cite{wannerer12}} for general Minkowski valuations (see also \textbf{\cite{SchuWan14}}).

\vspace{0.1cm}

In order to state the crucial integral representation of smooth translation invariant
and $\mathrm{SO}(n)$ equivariant Minkowski valuations we have to briefly recall
the convolution between functions and measures on $S^{n-1}$.
First note that since $\mathrm{SO}(n)$ is a compact Lie group, the convolution $\sigma \ast \tau$
of signed measures $\sigma, \tau$ on $\mathrm{SO}(n)$ can be defined by
\[\int_{\mathrm{SO}(n)}\!\!\! f(\vartheta)\, d(\sigma \ast \tau)(\vartheta)=\int_{\mathrm{SO}(n)}\! \int_{\mathrm{SO}(n)}\!\!\! f(\eta \theta)\,d\sigma(\eta)\,d\tau(\theta), \qquad f \in C(\mathrm{SO}(n)).   \]

Using the identification of the sphere $S^{n-1}$ with the homogeneous space $\mathrm{SO}(n)/\mathrm{SO}(n-1)$
leads to a one-to-one correspondence of
$C(S^{n-1})$ and $\mathcal{M}(S^{n-1})$ with right
$\mathrm{SO}(n-1)$ invariant functions and measures on $\mathrm{SO}(n)$ (cf.\
\textbf{\cite{Schu09}} for more information). Using this
correspondence, the convolution of measures on $\mathrm{SO}(n)$ induces a convolution product
on $\mathcal{M}(S^{n-1})$.

\pagebreak

For spherical convolution zonal functions and measures on
$S^{n-1}$ play an essential role. We denote the set of continuous zonal functions on $S^{n-1}$ by
$C(S^{n-1},\bar{e})$. For $\sigma \in \mathcal{M}(S^{n-1})$, $f
\in C(S^{n-1},\bar{e})$, and $\eta \in \mathrm{SO}(n)$, it is
easy to check that
\begin{equation} \label{zonalconv}
(\sigma \ast f)(\eta \bar{e})
=\int_{S^{n-1}}f(\eta^{-1} u)\,d\sigma(u).
\end{equation}
Note that, by (\ref{zonalconv}), we have, for every $\vartheta \in \mathrm{SO}(n)$, that
\begin{equation*}
(\vartheta \sigma) \ast f = \vartheta(\sigma \ast f),
\end{equation*}
where $\vartheta\sigma$ is the image measure of $\sigma$ under the rotation $\vartheta \in \mathrm{SO}(n)$. Moreover, from the obvious identification of
zonal functions on $S^{n-1}$ with functions on $[-1,1]$, (\ref{zonalconv}), and the Funk--Hecke Theorem, it follows that there
are $a_k^n[f] \in \mathbb{R}$ such that the spherical harmonic expansion
of $\sigma \ast f \in C(S^{n-1})$ is given by
\[\sigma \ast f \sim \sum_{k=0}^{\infty} a_k^n[f]\,\pi_k\sigma.  \]
Hence, convolution from the right induces a multiplier transformation. It is also not difficult to check from (\ref{zonalconv}) that the convolution of
zonal functions and measures is Abelian.

Another property of spherical convolution which is going to be critical for us is the fact that the convolution is selfadjoint, in particular, we have
for all $\sigma, \tau \in \mathcal{M}(S^{n-1})$ and every $f \in C(S^{n-1},\bar{e})$,
\begin{equation} \label{convselfad}
\int_{S^{n-1}} (\sigma \ast f)(u)\,d\tau(u) = \int_{S^{n-1}} (\tau \ast f)(u)\,d\sigma(u).
\end{equation}

We are now in a position to state a recent Hadwiger type theorem for smooth Minkowski valuations which is
the key to the proof of Theorem \ref{thm3}.

\begin{theorem} \label{repschuwan} {\bf (\!\!\cite{wannerer12, SchuWan14})} If $\Phi_j \in \mathbf{MVal}_j^{\infty}$, $j \in \{1, \ldots, n - 1\}$, then there exists a
unique $f \in C_{\mathrm{o}}^{\infty}(S^{n-1},\bar{e})$, called the \emph{generating function} of $\Phi_j$, such that for every $K \in \mathcal{K}^n$,
\begin{equation}
h(\Phi_jK,\cdot) = S_j(K,\cdot) \ast f.
\end{equation}
\end{theorem}

\pagebreak

\noindent {\bf Examples}:
\begin{enumerate}
\item[(a)] Kiderlen \textbf{\cite{kiderlen05}} proved (in a slightly different form) that if $\Phi_1 \in \mathbf{MVal}_1^{\infty}$, then there exists
a unique $g \in C_{\mathrm{o}}^{\infty}(S^{n-1},\bar{e})$ such that for every $K \in \mathcal{K}^n$,
\begin{equation} \label{kidrep}
h(\Phi_1 K,\cdot) = h(K,\cdot) \ast g.
\end{equation}
In order see how (\ref{kidrep}) is related to Theorem \ref{repschuwan},
we use (\ref{boxhks1}) and the fact that $\Box_n: C_{\mathrm{o}}^{\infty}(S^{n-1}) \rightarrow C_{\mathrm{o}}^{\infty}(S^{n-1})$ is a bijective multiplier transformation,
to obtain a function $f \in C^{\infty}(S^{n-1},\bar{e})$ with $\Box_n f = g$ and conclude that
\[h(\Phi_1 K,\cdot) = h(K,\cdot) \ast g = h(K,\cdot) \ast  \Box_n f = \Box_n h(K,\cdot) \ast f = S_1(K,\cdot) \ast f. \]

\item[(b)] The case $j = n - 1$ of Theorem \ref{repschuwan} was first proved (in a more general form) in \textbf{\cite{Schu06a}}.
Moreover, it was also shown there that $\Phi_{n-1} \in \mathbf{MVal}_{n-1}^{\infty}$ is \emph{even}
if and only if there exists an $o$-symmetric smooth convex body of revolution $L \in \mathcal{K}^n$ such that for every $K \in \mathcal{K}^n$,
\begin{equation*} \label{schurep}
h(\Phi_{n-1} K,\cdot) = S_{n-1}(K,\cdot) \ast h(L,\cdot).
\end{equation*}

\item[(c)] For $i \in \{1, \ldots, n - 1\}$, the support function of the projection body map of order $i$,  $\Pi_i \in \mathbf{MVal}_i$, is given by
\[h(\Pi_iK,u) = V_i(K|u^{\bot}) = \frac{1}{2} \int_{S^{n-1}} |u\cdot v|\,dS_i(K,v), \qquad u \in S^{n-1}.   \]
Note that $\Pi_i$ is continuous but {\it not} smooth. Its (merely) continuous generating function
is given by $f(u) = \frac{1}{2}|u \cdot \bar{e}|$, $u \in S^{n-1}$.

\item[(d)] For $i \in \{2, \ldots, n\}$, the (normalized) mean section operator of order $i$, $\mathrm{M}_i \in \mathbf{MVal}_{n+1-i}$,
was first defined in \textbf{\cite{GoodeyWeil1992}} by
\[h(\mathrm{M}_iK,\cdot) = \int_{\mathrm{AGr}_{i,n}} h(\mathrm{J}(K \cap E),\cdot)\,dE.  \]
Here, $\mathrm{J} \in \mathbf{MVal}_1$ is defined by $\mathrm{J}K =
K - s(K)$, where $s: \mathcal{K}^n \rightarrow \mathbb{R}^n$ \linebreak is
the Steiner point map. Recently, Goodey and Weil \textbf{\cite{GoodeyWeil2014}} proved that the
generating functions of the mean section operators are up to normalization the zonal functions $\breve{g}_i \in L_1(S^{n-1},\bar{e})$ determined by C.~Berg's functions $g_i$ on $[-1,1]$.
More precisely,
\begin{equation} \label{genfctmi}
h(\mathrm{M}_iK,\cdot) = p_{n,i}\,S_{n+1-i}(K,\cdot) \ast \breve{g}_i,
\end{equation}
with constants $p_{n,i}$ which were explicitly
determined in \textbf{\cite{GoodeyWeil2014}}.
\end{enumerate}

The integration operator $\mathfrak{L}$ on translation invariant scalar valuations can be
extended to Minkowski valuations by using (\ref{hardlefint}):
\[h((\mathfrak{L}\Phi)(K),\cdot) = \int_{\mathrm{AGr}_{n-1,n}} h(\Phi(K \cap E),\cdot)\,dE.  \]

It was proved in \textbf{\cite{papschu12}} that also the
derivation operator $\Lambda$ can be extended to continuous translation
invariant Minkowski valuations:
\[h((\Lambda \Phi)(K),\cdot) = \left . \frac{d}{dt} \right |_{t=0} h(\Phi(K + tB),\cdot).  \]
Note that in this case it is not trivial that the right hand side actually defines the support function of a convex body; this was proved in \textbf{\cite{papschu12}}.

If $\Phi_j \in \mathbf{MVal}_j^{\infty}$, $1 \leq j \leq n - 1$, with associated real valued valuation $\nu_j \in \mathbf{Val}_j^{\infty}$,
then the associated real valued valuations of $\mathfrak{L}\Phi_j$ and $\Lambda \Phi_j$ are given by
$\mathfrak{L}\nu_j \in \mathbf{Val}_{j+1}^{\infty}$ and $\Lambda \nu_j \in \mathbf{Val}_{j+1}^{\infty}$, respectively.
In particular, we have  $\mathfrak{L}\Phi_j \in \mathbf{MVal}_{j+1}^{\infty}$ and $\Lambda \Phi_j \in \mathbf{MVal}_{j-1}^{\infty}$.

In view of Theorem \ref{repschuwan}, it is a natural problem to
determine the induced action of the $\mathrm{SO}(n)$ equivariant
operators $\Lambda$ and $\mathfrak{L}$ on the generating
functions of smooth Minkowski valuations. This was done in
\textbf{\cite{SchuWan13}} and is the content of the following
theorem.

\begin{theorem} \label{actionhardlef} {\bf (\!\!\cite{SchuWan13})}
Suppose that $\Phi_j \in \mathbf{MVal}_j^{\infty}$ and let $f \in
C_{\mathrm{o}}^{\infty}(S^{n-1},\bar{e})$ be the generating
function of $\Phi_j$.
\begin{enumerate}
\item[(a)] If $2 \leq j \leq n - 1$, then the generating function
of $\Lambda \Phi_j$ is given by $jf$.
\item[(b)] If $1 \leq j \leq n - 2$, then there exists a constant $c_{n,j} > 0$ such that the generating function
of $\mathfrak{L}\Phi_j$ is given by $c_{n,j}\,\Box_{n-j+1} f \ast
\breve{g}_{n-j}$.
\end{enumerate}
In particular, $\Lambda: \mathbf{MVal}_j^{\infty} \rightarrow
\mathbf{MVal}_{j-1}^{\infty}$ is injective for all $2 \leq j \leq
n - 1$ and $\mathfrak{L}: \mathbf{MVal}_j^{\infty} \rightarrow
\mathbf{MVal}_{j+1}^{\infty}$ is injective for all $1 \leq j \leq
n - 2$.
\end{theorem}

The constants $c_{n,j}$ from Theorem \ref{actionhardlef} (b) were
explicitly determined in \textbf{\cite{SchuWan13}}. We also note
that, by Theorem \ref{actionhardlef}, for $i > j$, the map
\[\Lambda^{j-i}: \mathbf{MVal}_{j,i}^{\infty} \rightarrow \mathbf{MVal}_i^{\infty}  \]
is well defined. From Theorem \ref{actionhardlef} (a) and
Examples (a) and (b) above, we can also deduce more information
about the classes $\mathbf{MVal}_{j,i}^{\infty}$.

\pagebreak

\begin{koro} \label{mvalijkor} $\phantom{a}$
\begin{enumerate}
\item[(a)] Suppose that $1 \leq i,j \leq n - 1$, $\Phi_j \in
\mathbf{MVal}_j^\infty$, and let $f \in
C_{\mathrm{o}}^{\infty}(S^{n-1},\bar{e})$ be the generating
function of $\Phi_j$. Then $\Phi_j \in \mathbf{MVal}_{j,i}^\infty$
if and only if $S_i(K,\cdot) \ast f$ is a support function for
every $K \in \mathcal{K}^n$.
\item[(b)] $\mathbf{MVal}_{1,n-1}^\infty \varsubsetneq
\mathbf{MVal}_1^{\infty}$.
\end{enumerate}
\end{koro}
{\it Proof.} Statement (a) is a direct consequence of the
definition of $\mathbf{MVal}_{j,i}^\infty$ \linebreak and Theorem
\ref{actionhardlef} (a).

In order to prove (b) let $\Phi_1 \in \mathbf{MVal}_1^{\infty}$
be even and let $f \in C_{\mathrm{o}}^{\infty}(S^{n-1},\bar{e})$
be the generating function of $\Phi_j$. Then, by (a) and example
(b) from above, $\Phi_1 \in \mathbf{MVal}_{1,n-1}^\infty$ if and
only if $f=h(L,\cdot)$ for some $o$-symmetric smooth convex body
of revolution $L \in \mathcal{K}^n$. Thus, we have
\[h(\Phi_1 K,\cdot) = S_1(K,\cdot) \ast h(L,\cdot) = \Box_n h(K,\cdot) \ast h(L,\cdot) = h(K,\cdot) \ast  s_1(L,\cdot), \]
where $s_1(L,\cdot)$ is the smooth density of $S_1(L,\cdot)$. It
was proved by Kiderlen \textbf{\cite{kiderlen05}} that for any
(even) non-negative $g \in
C_{\mathrm{o}}^{\infty}(S^{n-1},\bar{e})$, (\ref{kidrep}) defines
an (even) Minkowski valuation in $\mathbf{MVal}_1^{\infty}$. Since
area measures of order 1 are nowhere dense in
$\mathcal{M}_{\mathrm{o}}$, this proves the claim. \hfill
$\blacksquare$

\vspace{0.3cm}

Note that, by Corollary \ref{mvalijkor} (b), in general
$\mathbf{MVal}_{j,i}^\infty \varsubsetneq
\mathbf{MVal}_j^{\infty}$ for $i > j$. Explicit examples of
Minkowski valuations in $\mathbf{MVal}_j$ with generating functions
which do not generate a Minkowski valuation in
$\mathbf{MVal}_{n-1}$ are provided by the mean section operators.
This follows from (\ref{genfctmi}) and the case $i = n - 1$ of
Theorem \ref{repschuwan} for continuous Minkowski valuations
established in \textbf{\cite{Schu06a}}, where it was proved that
$\Phi_{n-1} \in \mathbf{MVal}_{n-1}$ is generated by a
\emph{continuous} function $f \in
C_{\mathrm{o}}(S^{n-1},\bar{e})$. However, C.~Berg's functions
$g_i$ are not continuous on $[-1,1]$ for $i \geq 5$.

\vspace{0.1cm}

We end this section with another remark concerning Corollary
\ref{mvalijkor} (a): Generating functions or earlier versions of
Theorem \ref{repschuwan}, respectively, were the critical tool
used in the proofs of the first Brunn--Minkowski type inequalities
for Minkowski valuations. In the next section we will see that
the Hard Lefschetz operators on Minkowski valuations (which were
introduced only recently) and Theorem \ref{repschuwan} both
naturally lead to the same classes $\mathbf{MVal}_{j,i}$ for which
we can establish such inequalities.

\pagebreak

\centerline{\large{\bf{ \setcounter{abschnitt}{6}
\arabic{abschnitt}. Proofs of main results}}}

\reseteqn \alpheqn \setcounter{theorem}{0}

\vspace{0.6cm}

After these preparations, we are now in a position to complete
the proofs of Theorems \ref{thm3}, \ref{thm4}, and
\ref{mainorlicz}. We begin with Theorem \ref{thm3} which we first
restate.

\begin{theorem} \label{thm51} Let $\Phi_j \in
\mathbf{MVal}_j^{\infty}$, $2 \leq j \leq n - 1$. For $j + 2 \leq
i \leq n$ and every $L \in \mathcal{K}^n$, there exists
$\gamma_{i,j}(L,\cdot) \in \mathbf{Val}_1^{-\infty}$ such that
\begin{equation*}
W_{n-i}(K,\Phi_jL)=\gamma_{i,j}(L,(\mathfrak{L}^{i-j-1}\Phi_j)(K))
\end{equation*}
for every $K \in \mathcal{K}^n$. Moreover, if $\Phi_j \in
\mathbf{MVal}_{j,i-1}^{\infty}$, then
\[\gamma_{i,j}(L,(\mathfrak{L}^{i-j-1}\Phi_j)(K)) = \frac{(i-1)!}{j!} W_{n-1-j}(L,(\Lambda^{j+1-i}\Phi_j)(K)).   \]
\end{theorem}
{\it Proof.} First define an isomorphism $\Theta_j:
C_{\mathrm{o}}^{\infty}(S^{n-1}) \rightarrow
C_{\mathrm{o}}^{\infty}(S^{n-1})$ by
\[\Theta_jf = c_{n,j\,} \Box_{n-j+1} f \ast \breve{g}_{n-j} = c_{n,j\,} f \ast \Box_{n-j+1} \breve{g}_{n-j},  \]
where the constant $c_{n,j} > 0$ is as in Theorem
\ref{actionhardlef} (b). Here, the second equality follows from
the fact that multiplier transformations commute and $\Box_{n-j+1}
\breve{g}_{n-j}$ is to be understood in the sense of distributions,
where we use the canonical extension of the selfadjoint operator
$\Box_{n-j+1}$ to $C_{\mathrm{o}}^{-\infty}(S^{n-1})$.

Let $\tau_{\bar{e}} = \delta_{\bar{e}}- \pi_1\delta_{\bar{e}} \in
\mathcal{M}_{\mathrm{o}}(S^{n-1})$, where $\delta_{\bar{e}}$ is
the Dirac measure supported in $\bar{e} \in S^{n-1}$. Then, by
(\ref{zonalconv}), $f \ast \tau_{\bar{e}} = f$ for every $f \in
C_{\mathrm{o}}^{\infty}(S^{n-1})$. Now since $\Box_k \breve{g}_k
= \tau_{\bar{e}}$, it follows from Theorem \ref{actionhardlef}
(b) that if $f \in C_{\mathrm{o}}^{\infty}(S^{n-1},\bar{e})$ is
the generating function of $\Phi_j$, then
$\mathfrak{L}^{i-j-1}\Phi_j \in \mathbf{MVal}_{i-1}^{\infty}$ is
generated by
\begin{equation} \label{thetasucc}
\Theta_{i-2} \Theta_{i-1} \cdots \Theta_{j+1}\Theta_j f =
q_{n,i,j\,} \Box_{n-j+1}f \ast \breve{g}_{n-i+2},
\end{equation}
where $q_{n,i,j} = \prod_{m=j}^{i-2}c_{n,m} > 0$. Note that the
inverse of the isomorphism (\ref{thetasucc}) is, for $f \in
C_{\mathrm{o}}^{\infty}(S^{n-1})$, given by
\[q_{n,i,j\,}^{-1} \Box_{n-i+2} f \ast \breve{g}_{n-j+1}. \]

For every $L \in \mathcal{K}^n$ we define a distribution
$T_{i,j}(L) \in C^{-\infty}_{\mathrm{o}}(S^{n-1})$ by
\[\langle f, T_{i,j}(L) \rangle = q_{n,i,j\,}^{-1} \int_{S^{n-1}} (\Box_{n-i+2} f \ast \breve{g}_{n-j+1})(u)\,dS_j(L,u) \]
for $f \in C^{-\infty}_{\mathrm{o}}(S^{n-1})$. Let
$\gamma_{i,j}(L,\cdot) \in \mathbf{Val}_1^{-\infty}$ be the
generalized valuation corresponding to $T_{i,j}(L)$ determined by
Theorem \ref{genvalchar}.

\pagebreak

Since $\mathfrak{L}^{i-j-1}\Phi_j$ is smooth, it follows that
$h((\mathfrak{L}^{i-j-1}\Phi_j)(K),\cdot)$ is smooth for every $K
\in \mathcal{K}^n$. Hence, we can evaluate
$\gamma_{i,j}(L,\cdot)$ on $(\mathfrak{L}^{i-j-1}\Phi_j)(K)$.
Using that
\[h((\mathfrak{L}^{i-j-1}\Phi_j)(K),\cdot) =  q_{n,i,j\,} S_{i-1}(K,\cdot) \ast ( \Box_{n-j+1}f \ast \breve{g}_{n-i+2}), \]
we obtain
\begin{eqnarray*}
\gamma_{i,j}(L,(\mathfrak{L}^{i-j-1}\Phi_j)(K)) & = & \langle h((\mathfrak{L}^{i-j-1}\Phi_j)(K),\cdot), T_{i,j}(L) \rangle \\
 & = & \int_{S^{n-1}} (S_{i-1}(K,\cdot) \ast  f)(u)\,dS_j(L,u).
\end{eqnarray*}
Now on one hand it follows from (\ref{convselfad}), that
\[\gamma_{i,j}(L,(\mathfrak{L}^{i-j-1}\Phi_j)(K)) =
\int_{S^{n-1}}\!\!\! (S_j(L,\cdot)\, {\ast} f)(u)\,dS_{i-1}(K,u)
=  W_{n-i}(K,\Phi_jL).\] On the other hand, if $\Phi_j \in
\mathbf{MVal}_{j,i-1}^{\infty}$, then, by Theorem
\ref{actionhardlef} (a),
\[ S_{i-1}(K,\cdot) \ast  f =  \frac{(i-1)!}{j!} h((\Lambda^{j+1-i}\Phi_j)(K),\cdot) \]
and, thus,
\[\gamma_{i,j}(L,(\mathfrak{L}^{i-j-1}\Phi_j)(K)) = \frac{(i-1)!}{j!} W_{n-1-j}(L,(\Lambda^{j+1-i}\Phi_j)(K)).  \]
which completes the proof. \hfill $\blacksquare$

\vspace{0.4cm}

Note that, by Theorem \ref{multipliers} and Lemma
\ref{distharmonic}, $\Box_{n-i+2} \breve{g}_{n-j+1} \in
C^{-\infty}_{\mathrm{o}}(S^{n-1})$ and that if $f \in
C_{\mathrm{o}}(S^{n-1})$, then also
\begin{equation} \label{actualdist}
f \ast \Box_{n-i+2} \breve{g}_{n-j+1} = \Box_{n-i+2} f \ast
\breve{g}_{n-j+1} \in C^{-\infty}_{\mathrm{o}}(S^{n-1}).
\end{equation}
However, in general (\ref{actualdist}) does \emph{not} define a
continuous function on $S^{n-1}$ if $f$ is merely continuous.

\vspace{0.1cm}

Next, we note that using Theorem \ref{repschuwan} we can also give
a new and short proof of Theorem \ref{thm2}: If $\Phi_j \in
\mathbf{MVal}_j^{\infty}$, $2 \leq j \leq n - 1$, has generating
function $f \in C^{\infty}_{\mathrm{o}}(S^{n-1})$ and $1 \leq i
\leq j + 1$, then, by (\ref{convselfad}) and Theorem
\ref{actionhardlef} (a),
\begin{eqnarray*}
W_{n-i}(K,\Phi_jL) & = & \int_{S^{n-1}} (S_{i-1}(K,u) \ast f)(u)\,dS_j(L,\cdot) \\
 & = &  \frac{(i-1)!}{j!}\,W_{n-j-1}(L,(\Lambda^{j+1-i}\Phi_j)(K)).
\end{eqnarray*}
for every $K, L \in \mathcal{K}^n$.

\pagebreak

Putting together Theorem \ref{thm2} and Theorem \ref{thm51} we
obtain the following.

\begin{koro} For $1 \leq i \leq n$, $2 \leq j \leq n - 1$, and $\Phi_j
\in \mathbf{MVal}_{j,i-1}^{\infty}$, we have
\begin{equation} \label{switch}
W_{n-i}(K,\Phi_jL)=\frac{(i-1)!}{j!}\,W_{n-1-j}(L,(\Lambda^{j+1-i}\Phi_j)(K))
\end{equation}
for every $K, L \in \mathcal{K}^n$.
\end{koro}

For the proof of Theorem \ref{mainorlicz} and in order to establish the equality cases in Theorem \ref{thm4}, we need the following monotonicity property
of Minkowski valuations:

\begin{lem} \label{monolem}
Suppose that $1 \leq i \leq n$, $2\leq j\leq n-1$, and let $\Phi_j
\in \mathbf{MVal}_{j,i-1}$ be non-trivial. If $K,L \in
\mathcal{K}^n$ have non-empty interiors, then $K \subseteq L$
implies that
\begin{equation}\label{desire1}
W_{n-i}(\Phi_jK)\leq W_{n-i}(\Phi_jL)
\end{equation}
with equality if and only if $K=L$. In particular, $W_{n-i}(\Phi_jK)>0$ for every $K \in \mathcal{K}^n$ with non-empty interior.
\end{lem}
{\it Proof.} We first assume that $i \geq 2$ and that $\Phi_j$ is
smooth. By (\ref{switch}) and the monotonicity of mixed volumes,
we have for every $Q \in \mathcal{K}^n$,
\begin{eqnarray*}
W_{n-i}(Q,\Phi_jL) & = & \frac{(i-1)!}{j!}W_{n-1-j}(L,(\Lambda^{j+1-i}\Phi_j)(Q))\\
& \geq & \frac{(i-1)!}{j!}W_{n-1-j}(K,(\Lambda^{j+1-i}\Phi_j)(Q))
=  W_{n-i}(Q,\Phi_jK).
\end{eqnarray*}
Thus, taking $Q = \Phi_j L$ and using inequality (\ref{genmink}),
yields
\[W_{n-i}(\Phi_jL)^i \geq W_{n-i}(\Phi_j L,\Phi_jK)^i \geq W_{n-i}(\Phi_j L)^{i-1}W_{n-i}(\Phi_jK)\]
which implies the desired inequality (\ref{desire1}). If $\Phi_j
\in \mathbf{MVal}_{j,i-1}$ is not smooth, (\ref{desire1}) follows
by approximation.

In order to establish the equality conditions first note that, by
the $\mathrm{SO}(n)$ equivariance of $\Phi_j$, the convex body
$\Phi_jB$ must be an $o$-symmetric ball. Moreover, from Lemma
\ref{basicMinkVal} (b), it follows that $\Phi_j B = r_{\Phi_j}B$,
where $r_{\Phi_j}
> 0$. Thus, since $K$ and $L$ have non-empty interiors, we conclude
from (\ref{desire1}) that $W_{n-i}(\Phi_j K), W_{n-i}(\Phi_j L)>0$
or, equivalently, that $\Phi_j K$ and $\Phi_j L$ have dimension
at least $i$.

Assume now that equality holds in (\ref{desire1}). Then, by the
equality conditions of (\ref{genmink}) and Lemma
\ref{basicMinkVal} (a), there exists an $\alpha
> 0$ such that $\Phi_j K = \alpha\, \Phi_j L$.
It follows from equality in (\ref{desire1}) that $\alpha = 1$.
Thus, by Lemma \ref{basicMinkVal} (b), we~have
\begin{equation}\label{eq:glkl}
W_{n-j}(K) = r_{\Phi_j}^{-1} W_{n-1}(\Phi_j K) = r_{\Phi_j}^{-1} W_{n-1}(\Phi_j L) = W_{n-j}(L).
\end{equation}
Using again the monotonicity of mixed volumes and
\eqref{genmink}, we obtain
\[W_{n-j}(L)^{j}  = W_{n-j}(L,L)^{j} \geq W_{n-j}(L,K)^{j} \geq  W_{n-j}(L)^{j-1}W_{n-j}(K).\]
From (\ref{eq:glkl}) and the equality conditions of inequality
(\ref{genmink}), we conclude that $K$ is a translate of $L$. But since $K\subseteq L$, we must have $K=L$.

Inequality (\ref{desire1}) for $i = 1$ follows directly from Lemma
\ref{basicMinkVal} (b) and the monotonicity of quermassintegrals.
If equality holds in (\ref{desire1}) for $i = 1$, then we have
(\ref{eq:glkl}) and therefore, as before, obtain that $K = L$.
\hfill $\blacksquare$

\vspace{0.4cm}

In contrast to Lemma \ref{monolem}, we note that not every
Minkowski valuation $\Phi_j \in \mathbf{MVal}_{j,i-1}$ is monotone
with respect to set inclusion (cf. \textbf{\cite{kiderlen05}}).
However, all known examples of Minkowski valuations $\Phi_j \in
\mathbf{MVal}_j$, $1 \leq j \leq n - 1$ are \emph{weakly
monotone}, that is, for every pair of convex bodies $K, L \in
\mathcal{K}^n$ such that $K \subseteq L$, there exists a vector
$x(K,L) \in \mathbb{R}^n$ such that
\[\Phi_j K \subseteq \Phi_j L + x(K,L).  \]
It is an open problem whether all translation invariant and
$\mathrm{SO}(n)$ equivariant Minkowski valuations are weakly
monotone. Using arguments as in the proof of Lemma \ref{monolem},
we can show the following.

\begin{prop} Suppose that $2 \leq j \leq n - 1$. If $\Phi_j
\in \mathbf{MVal}_{j,n-1}$, then $\Phi_j$ is weakly monotone.
\end{prop}
{\it Proof.} Without loss of generality we may assume that
$\Phi_j$ is smooth. If $K, L \in \mathcal{K}^n$ such that $K
\subseteq L$, then, as in Lemma \ref{monolem}, it follows from
(\ref{switch}) and the monotonicity of mixed volumes that for
every $Q \in \mathcal{K}^n$,
\begin{eqnarray*}
W_0(Q,\Phi_jL) & = & \frac{(n-1)!}{j!}W_{n-1-j}(L,(\Lambda^{j+1-n}\Phi_j)(Q))\\
& \geq & \frac{(n-1)!}{j!}W_{n-1-j}(K,(\Lambda^{j+1-n}\Phi_j)(Q))
=  W_0(Q,\Phi_jK).
\end{eqnarray*}
But it is well known (cf.\ \textbf{\cite[\textnormal{Corollary
4.3}]{Schu06a}}) that $W_0(Q,\Phi_jK) \leq W_0(Q,\Phi_jL)$ for
every $Q \in \mathcal{K}^n$ implies $\Phi_jK \subseteq \Phi_j L +
x$ for some $x \in \mathbb{R}^n$. \hfill $\blacksquare$

\vspace{0.4cm}

\pagebreak

We return now to the proof of Theorem \ref{thm4} which we also restate.

\begin{theorem} \label{th:mainth}
Let $1 \leq i \leq n$ and let $\Phi_j \in \mathbf{MVal}_{j,i-1}$, $2 \leq j \leq n - 1$, be non-trivial. If $K, L \in \mathcal{K}^n$ have non-empty interiors,
then for all $\lambda \in (0,1)$,
\begin{equation} \label{mainineq17}
W_{n-i}(\Phi_j((1-\lambda)K+\lambda L)) \geq
W_{n-i}(\Phi_jK)^{1-\lambda}W_{n-i}(\Phi_jL)^\lambda,
\end{equation}
with equality if and only if $K$ and $L$ are translates of each
other.
\end{theorem}
{\it Proof.} First we assume again that $i \geq 2$ and that $\Phi_j$ is
smooth. We also use the abbreviations $K_\lambda = (1-\lambda) K + \lambda L$ and $Q = \Phi_j K_\lambda$.
Then, by (\ref{switch}),
\[W_{n-i}(\Phi_j K_\lambda) = W_{n-i}(Q,\Phi_j K_\lambda) = \frac{(i-1)!}{j!} W_{n-1-j}(K_\lambda,(\Lambda^{j+1-i}\Phi_j)(Q)).\]
From an application of inequality (\ref{mostgenbm}), we therefore obtain
\[W_{n-i}(\Phi_j K_\lambda)\! \geq\!\! \frac{(i{-}1)!}{j!}  W_{n\mbox{-}1\mbox{-}j}(K,\!(\Lambda^{j+1-i}\Phi_j)(Q))^{1-\lambda}W_{n\mbox{-}1\mbox{-}j}(L,\!(\Lambda^{j+1-i}\Phi_j)(Q))^{\lambda}.\]
Thus, using (\ref{switch}) again followed by (\ref{genmink}), we obtain
\begin{eqnarray*}
W_{n-i}(\Phi_jK_\lambda)^i & \geq & W_{n-i}(Q,\Phi_jK)^{i(1-\lambda)}W_{n-i}(Q,\Phi_jL)^{i\lambda} \\
 & = & W_{n-i}(Q)^{i-1}W_{n-i}(\Phi_jK)^{1-\lambda}W_{n-i}(\Phi_jL)^{\lambda}.
\end{eqnarray*}
Since $Q=\Phi_j K_\lambda$, this is the desired inequality (\ref{mainineq17}). If $\Phi_j
\in \mathbf{MVal}_{j,i-1}$ is not smooth, (\ref{mainineq17}) follows now
by approximation.

In order to establish the equality conditions first note that by Lemma \ref{monolem},
$\Phi_jK$, $\Phi_jL$, and $\Phi_jK_{\lambda}$ all have dimension at least $i$. Therefore,
the equality conditions of inequality (\ref{genmink}) imply that $\Phi_jK$ is homothetic to $\Phi_jK_\lambda$, which is in turn homothetic to $\Phi_jL$.
In fact, by Lemma \ref{basicMinkVal} (a), they have to be dilates of one another, that is, there exist
$t_1, t_2 > 0$ such that
\[t_1 \Phi_j K = \Phi_j K_\lambda =t_2 \Phi_j L,\]
where $1=t_1^{1-\lambda}t_2^\lambda$, by the equality in (\ref{mainineq17}). Moreover, an application of Lemma \ref{basicMinkVal} (b)
yields $t_1W_{n-j}(K)=W_{n-j}(K_\lambda)= t_2 W_{n-j}(L)$.
Consequently, we have
\[W_{n-j}(K_\lambda) = W_{n-j}(K)^{1-\lambda} W_{n-j}(L)^{\lambda}. \]
By the equality conditions of inequality (\ref{quermassbm}), this is possible only if $K$ and $L$ are translates. This completes
the proof for $i \geq 2$. If $i = 1$, then the statement is an immediate consequence of Lemma \ref{basicMinkVal} (b) and (\ref{quermassbm}). \hfill $\blacksquare$

\pagebreak

It remains to complete the proof of Theorem \ref{mainorlicz}.

\begin{theorem} \label{th:orlicz} Let $\varphi \in \Theta_1$, $1 \leq i \leq n$, and let $\Phi_j \in
\mathbf{MVal}_{j,i-1}$, $2 \leq j \leq n - 1$, be non-trivial. If
$K, L \in \mathcal{K}^n$ contain the origin, then for all
$\lambda \in (0,1)$,
\begin{equation} \label{yoda1717}
W_{n-i}(\Phi_j(K+_{\varphi,\lambda} L)) \geq W_{n-i}(\Phi_jK)^{1-\lambda}W_{n-i}(\Phi_jL)^\lambda.
\end{equation}
When $\varphi$ is strictly convex and $K$ and $L$ have non-empty
interiors, equality holds if and only if $K = L$.
\end{theorem}
{\it Proof.} The case $i = 1$ follows directly from Lemma \ref{orliczinclusion}, the monotonicity of mixed volumes, and inequality (\ref{quermassbm}).
Similarly, inequality (\ref{yoda1717}) for $i \geq 2$ follows from Lemma \ref{orliczinclusion}, Lemma \ref{monolem}, and Theorem \ref{th:mainth}.

In order to establish the equality conditions for $i \geq 2$, let
$\varphi$ be strictly convex and let $K$ and $L$ have non-empty interiors.
It follows from the equality conditions of Lemma \ref{monolem} that
\begin{equation}\label{eq:cbe}
K +_{\varphi,\lambda}  L=(1-\lambda)K+ \lambda  L.
\end{equation}
We want to show that this is possible only if $K = L$ or, equivalently, if $h(K,u) = h(L,u)$ for all
$u\in S^{n-1}$. If $h(K,u)=h(L,u)=0$, then there is nothing to prove. Therefore, we may assume that
$h(K +_{\varphi,\lambda}  L,u )>0$.
Now from the definition of the Orlicz convex combination, (\ref{eq:cbe}), together with the convexity of $\varphi$ and our assumption that $\varphi(1)=1$, we obtain
\[\varphi\left(\frac{(1-\lambda)h(K,u)+\lambda h(L,u)}{h(K+_{\varphi,\lambda}  L,u )}\right) = 1.\]
But since we have assumed that $\varphi$ is strictly convex, this implies that $h(K,u)=h(L,u)$.
\hfill $\blacksquare$

\vspace{0.4cm}

Like the classical inequality (\ref{quermassbm}), Theorem \ref{th:mainth} as well as Theorem \ref{th:orlicz} in case of a homogeneous addition are equivalent to corresponding
additive versions. Here we state one such additive version for $L_p$ Minkowski addition.

\begin{koro} \label{additivekoro}
Let $p > 1$, $1 \leq i \leq n$, and let $\Phi_j \in
\mathbf{MVal}_{j,i-1}$, $2 \leq j \leq n - 1$, be non-trivial. If
$K, L \in \mathcal{K}^n$ contain the origin in their interiors, then
\begin{equation*}
V_i(\Phi_j((1-\lambda)\cdot K +_p \lambda\cdot L))^{\frac{p}{ij}} \geq
(1-\lambda) V_i(\Phi_jK)^{\frac{p}{ij}} + \lambda V_i(\Phi_jL)^{\frac{p}{ij}},
\end{equation*}
with equality if and only if $K$ and $L$ are dilates of each
other.
\end{koro}

We finally remark that the special case $j = n - 1$ of Corollary \ref{additivekoro} was recently obtained by Wang \textbf{\cite{Wang2013}}.

\pagebreak

\vspace{1cm}

\centerline{\large{\bf{\setcounter{abschnitt}{7} Appendix}}}

\vspace{0.2cm}

\centerline{\phantom{a}by Semyon Alesker}

\renewcommand{\theequation}{A.\arabic{equation}}{\setcounter{equation}{0}
\setcounter{theorem}{0}
\renewcommand{\thesection}{\Alph{section}}
\setcounter{section}{1}

\vspace{0.6cm}

The purpose of this appendix is to provide a proof of Theorem \ref{genvalchar}. To this end we first show that all valuations in
$\mathbf{Val}_1^{\infty}$ and $\mathbf{Val}_{n-1}^{\infty}$ are
of the form (\ref{lesnope1}). In order to prove this,
we want to apply the Irreducibility Theorem as well as a deep
result from representation theory by Casselmann--Wallach
\textbf{\cite{cassel89}}. Therefore, we need to rewrite the valuations
$\nu_{1,f}$ and $\nu_{n-1,f}$ in $\mathrm{GL}(n)$ invariant terms
without referring to a Euclidean structure.

Recall that a \emph{line bundle} over a smooth manifold $M$ consists of a
smooth manifold $E$ and a surjective smooth map $\pi: E
\rightarrow M$ satisfying the following:
\begin{itemize}
\item For each $p \in M$, the fiber $E_p = \pi^{-1}(p)$ is a
1-dimensional vector space.
\item Every $p \in M$ has an open neighborhood $U$ in $M$ for
which there exists a diffeomorphism $\varrho: \pi^{-1}(U)
\rightarrow U \times \mathbb{R}$ such that, for each $q \in U$,
$\varrho(E_q) \subseteq \{q\} \times \mathbb{R}$ and
$\varrho|_{E_q}: E_q \rightarrow \{p\} \times \mathbb{R}$ is a
linear isomorphism.
\end{itemize}
For more information on line bundles, in particular, the definitions of the dual of a line bundle
and the tensor product of line bundles needed in the following,
see, e.g., \textbf{\cite[\textnormal{p.\ 4}]{wallach73}}.

A \emph{section} of a line bundle $\pi: E \rightarrow M$ is a
continuous map $h: M \rightarrow E$ such that $h(p) \in E_p$ for
every $p \in M$. We denote by $C(M,E)$ the vector space of all
sections of $E$ and by $C^{\infty}(M,E)$ the space of smooth
sections of $E$ endowed with the natural locally convex topology
which makes it a Fr\'echet space (see, e.g,
\textbf{\cite[\textnormal{Chapter 3}]{grosseretal}}). A sequence
of smooth sections converges in this topology if and only if in
local coordinates all the derivatives converge uniformly on compact
subsets.

Important examples of line bundles are density bundles of
manifolds (see,~e.g., \textbf{\cite[\textnormal{p.\
429}]{lee2nd}}). Recall that a \emph{density} on an
$n$-dimensional vector space $V$ is a function on the
$n$-fold product of $V$, $\delta: V \times \cdots \times V
\rightarrow \mathbb{R}$, such that if $A: V \rightarrow V$ is any
linear map, then
\[\delta(Av_1,\ldots,Av_n) = |\det A|\,\delta(v_1,\ldots,v_n).  \]
We denote by $\mathscr{D}(V)$ the vector space of all densities on $V$ and, as usual,
write $\Lambda^k(V)$ for the $k$th exterior power of $V$.

\begin{approp} \label{propdensity} {\bf ($\!\!$\cite[\textnormal{p.\ 428}]{lee2nd})} The vector space $\mathscr{D}(V)$
is $1$-dimensional and spanned by $|\omega|$ for any nonzero $\omega \in \Lambda^n(V^*)$.
\end{approp}

\pagebreak

The \emph{density bundle} $\pi: \mathscr{D}M \rightarrow M$ of a
smooth manifold $M$ is defined by
\[\mathscr{D}M = \coprod_{p \in M} \mathscr{D}(T_pM),  \]
where $\pi$ is the natural projection map taking each element of
$\mathscr{D}(T_pM)$ to $p$.

A \emph{density} on $M$ is a section of $\mathscr{D}M$. By
Proposition \ref{propdensity}, any nonvanishing $n$-form $\omega$
on $M$ determines a positive density $|\omega|$ on $M$. In fact,
if $\omega$ is a nonvanishing $n$-form on an open subset $U
\subseteq M$, then any density $\delta$ on $U$ is of the form
$\delta = f\,|\omega|$ for some continuous function $f$ on $U$.
\linebreak From this, it is now straightforward to define the
integral over $M$ of a compactly supported density on $M$.  We
refer to \textbf{\cite[\textnormal{p.\ 431ff}]{lee2nd}} for the
details.

If $M$ is a compact smooth manifold and $\pi: E \rightarrow M$ a
line bundle over $M$, then, using integration of densities on $M$,
one can define a canonical and nondegenerate bilinear pairing
\begin{equation} \label{pairing}
\langle \,\cdot\,,\,\cdot\, \rangle:  C(M,E) \times C(M,E^*
\otimes \mathscr{D}M) \rightarrow \mathbb{R}, \quad (f,g) \mapsto
\int_M [f,g].
\end{equation}
Here, $C(M,E^* \otimes \mathscr{D}M)$ is the space of sections of
the line bundle
\[E^* \otimes  \mathscr{D}M = \mathrm{Hom}(E,\mathscr{D}M), \]
whose fiber at $p \in M$ is the space of all linear maps
$E_p \rightarrow \mathscr{D}(T_pM)$, and
\[ [ \,\cdot\,,\,\cdot\, ]: C(M,E) \times C(M,E^* \otimes \mathscr{D}M) \rightarrow C(M,\mathscr{D}M)  \]
is pointwise just the evaluation map.

Now let $V$ be an $n$-dimensional vector space and denote by
$\mathbb{P}_+(V^*)$ the oriented projectivized cotangent bundle,
that is, the compact smooth manifold given by
\[\mathbb{P}_+(V^*) = (V^*\backslash \{0\})/\mathbb{R}^+.  \]
In the following we write $[\xi] := \mathrm{span}\, \xi$ for the
$1$-dimensional linear span of $\xi \in V^*\backslash \{0\}$ and
we use $[\xi]_+$ to denote the elements of $\mathbb{P}_+(V^*)$.
Note that if we choose a Euclidean structure on $V$, then we can
identify $V^*$ with $V$ and $\mathbb{P}_+(V^*)$ is diffeomorphic
to $S^{n-1}$. However, in contrast to $S^{n-1}$, we have a
natural $\mathrm{GL}(V)$ action on $\mathbb{P}_+(V^*)$ given by
\[A \cdot [\xi]_+ = [A\cdot\xi]_+, \qquad A \in \mathrm{GL}(V),  \]
where $A\cdot\xi \in V^*$ is defined by $(A\cdot\xi)(v) =
\xi(A^{-1}v)$ for $v \in V$. Also note that $\mathrm{GL}(V)$ acts
naturally on $\mathscr{D}(V)$ by
\[A\cdot \delta(v_1,\ldots,v_n) = \delta(A^{-1}v_1,\ldots,A^{-1}v_n), \qquad  A \in \mathrm{GL}(V), v_j \in V. \]

\pagebreak

\begin{apthm} \label{thma1} $\phantom{a}$
\begin{enumerate}
\item[(a)] The map $C^{\infty}_{\mathrm{o}}(S^{n-1}) \rightarrow
\mathbf{Val}_1^{\infty}$, $f \mapsto \nu_{1,f}$, where
\begin{equation} \label{nu1f}
\nu_{1,f}(K)= \int_{S^{n-1}} f(u)\,dS_1(K,u),
\end{equation}
is an isomorphism of Fr\'echet spaces.
\item[(b)] The map $C^{\infty}_{\mathrm{o}}(S^{n-1}) \rightarrow
\mathbf{Val}_{n-1}^{\infty}$, $f \mapsto \nu_{n-1,f}$, where
\begin{equation} \label{nun1f}
\nu_{n-1,f}(K)= \int_{S^{n-1}} f(u)\,dS_{n-1}(K,u),
\end{equation}
is an isomorphism of Fr\'echet spaces.
\end{enumerate}
\end{apthm}
{\it Proof.} First we note that, by the density properties of area
measures, both maps $f \mapsto \nu_{1,f}$ and $f \mapsto
\nu_{n-1,f}$ are injective. Moreover, it is not difficult to show
that they are also both continuous in the respective Fr\'echet
topologies.

In order to prove (a), we use (\ref{boxhks1}) and the fact that
$\Box_n$ is self-adjoint, to rewrite (\ref{nu1f}) to
\[\nu_{1,f}(K)= \int_{S^{n-1}} \Box_n f(u)\, h(K,u)\, du.  \]
Since $\Box_n: C^{\infty}_{\mathrm{o}}(S^{n-1}) \rightarrow
C^{\infty}_{\mathrm{o}}(S^{n-1})$ is an isomorphism, it suffices
to show \nolinebreak that
\begin{equation} \label{fuhku}
f \mapsto \left ( K \mapsto \int_{S^{n-1}} f(u)\, h(K,u)\, du \right)
\end{equation}
is an isomorphism between $C^{\infty}_{\mathrm{o}}(S^{n-1})$ and
$\mathbf{Val}_1^{\infty}$. To this end, recall that the support
function of a convex body $K \in \mathcal{K}^n$ is a
$1$-homogeneous function on $V^*$ and, thus, can be identified
with a section of a line bundle $E$ over $\mathbb{P}_+(V^*)$ whose
fiber over $[\xi]_{+} \in \mathbb{P}_+(V^*)$ is given by
$E_{[\xi]_+} = \mathscr{D}([\xi])$. To be more precise, we
identify $h(K,\cdot)$ with the section $\bar{h}(K,\cdot) \in
C(\mathbb{P}_+(V^*),E)$ defined by
\[\bar{h}(K,[\xi]_+)(c\xi) = |c|h(K,\xi), \qquad \xi \in V^*,\, c \in \mathbb{R}.   \]

Note that if we choose a Euclidean structure on $V$, then
$C(\mathbb{P}_+(V^*),E)$ can be identified with $C(S^{n-1})$. In
the same way, the Fr\'echet space of smooth sections
$C^{\infty}(\mathbb{P}_+(V^*),E^* \otimes
\mathscr{D}\,\mathbb{P}_+(V^*))$ is isomorphic to
$C^{\infty}(S^{n-1})$ and we let
$C^{\infty}_{\mathrm{o}}(\mathbb{P}_+(V^*),E^* \otimes
\mathscr{D}\,\mathbb{P}_+(V^*))$ denote the subspace isomorphic
to $C^{\infty}_{\mathrm{o}}(S^{n-1})$.

Using these identifications and the pairing defined in
(\ref{pairing}), the integral in (\ref{fuhku}) can be rewritten as
\begin{equation} \label{pairing1}
\int_{S^{n-1}} f(u)\, h(K,u)\, du = \int_{\mathbb{P}_+(V^*)}
[\bar{h}(K,\cdot),\bar{f}] = \langle \bar{h}(K,\cdot),
\bar{f}\rangle,\end{equation} where $\bar{f} \in
C^{\infty}_{\mathrm{o}}(\mathbb{P}_+(V^*),E^* \otimes
\mathscr{D}\,\mathbb{P}_+(V^*))$ denotes the section
corresponding to $f \in C^{\infty}_{\mathrm{o}}(S^{n-1})$.
Finally, note that the group $\mathrm{GL}(V)$ acts on the spaces
$C(\mathbb{P}_+(V^*),E)$ and
$C^{\infty}_{\mathrm{o}}(\mathbb{P}_+(V^*),E^* \otimes
\mathscr{D}\,\mathbb{P}_+(V^*))$ by left translation, that is,
\[(A\cdot \bar{f})([\xi]_+) = \bar{f}(A^{-1} \cdot [\xi]_+), \qquad A \in \mathrm{GL}(V),\,[\xi]_+ \in \mathbb{P}_+(V^*), \]
and the pairing (\ref{pairing1}) is invariant under these actions.
Thus, the map
\begin{equation} \label{inviso1}
C^{\infty}_{\mathrm{o}}(\mathbb{P}_+(V^*),E^* \otimes
\mathscr{D}\,\mathbb{P}_+(V^*)) \rightarrow
\mathbf{Val}_1^{\infty}, \quad \bar{f} \mapsto \langle
\bar{h}(K,\cdot),\bar{f} \rangle
\end{equation}
is $\mathrm{GL}(V)$ equivariant. Hence, its image is a
$\mathrm{GL}(V)$ invariant subspace of $\mathbf{Val}_1^{\infty}$
and, therefore, dense by the Irreducibility Theorem. However, by
the Casselmann--Wallach theorem \textbf{\cite{cassel89}}, this
image is also closed which proves \nolinebreak (a).

\vspace{0.1cm}

For the proof of (b), first note that if $P \in \mathcal{K}^n$ is
a polytope, then by the definition of $S_{n-1}(P,\cdot)$, we have
\begin{equation} \label{mcmullpoly}
\nu_{n-1,f}(P)= \sum_{F \in
\mathcal{F}_{n-1}(P)}\!\!\!f(u_F)\,\mathrm{vol}_{n-1}(F),
\end{equation}
where $\mathcal{F}_{n-1}(P)$ is the set of all facets of $P$ and
$u_F$ is the outer unit normal vector of the facet $F$.
Conversely, it is well known that if $f \in
C_{\mathrm{o}}(S^{n-1})$, then any function on polytopes in
$\mathbb{R}^n$ of the form (\ref{mcmullpoly}) has a unique
extension to a valuation in $\mathbf{Val}_{n-1}$ (see, e.g.,
\textbf{\cite[\textnormal{Chapter 6.4}]{schneider93}}). Moreover,
if $f \in C_{\mathrm{o}}^{\infty}(S^{n-1})$, then this valuation
is smooth and given by (\ref{nun1f}). In order to rewrite
$\nu_{n-1,f}$ in $\mathrm{GL}(V)$ invariant terms, it therefore
suffices to rewrite (\ref{mcmullpoly}).

To this end, let $\mathbb{P}_+^{\vee}(V)$ denote the compact
manifold of all cooriented $n-1$ dimensional subspaces in $V$.
(Recall that if $H$ is such a subspace, an orientation of $V/H$
is called a coorientation of $H$.) Note that there is a natural
diffeomorphism between $\mathbb{P}_+^{\vee}(V)$ and
$\mathbb{P}_+(V^*)$ which is equivariant under the action of
$\mathrm{GL}(V)$ on both manifolds. Let $S$ denote the line
bundle over $\mathbb{P}_+^{\vee}(V)$ whose fiber over $H \in
\mathbb{P}_+^{\vee}(V)$ is given by $S_H = \mathscr{D}(H)$. If we
choose a Euclidean structure,
$C^{\infty}(\mathbb{P}_+^{\vee}(V),S)$ is clearly isomorphic to
$C^{\infty}(S^{n-1})$ and we write again
$C^{\infty}_{\mathrm{o}}(\mathbb{P}_+^{\vee}(V),S)$ for the
subspace isomorphic to $C^{\infty}_{\mathrm{o}}(S^{n-1})$.

We now rewrite (\ref{mcmullpoly}) in the form
\begin{equation} \label{nun1inv1}
\nu_{n-1,\bar{f}}(P)= \sum_{F \in \mathcal{F}_{n-1}(P)}\! \int_{\hat{F}}
\bar{f}(\hat{F}),
\end{equation}
where $\bar{f} \in
C^{\infty}_{\mathrm{o}}(\mathbb{P}_+^{\vee}(V),S)$ is the section
corresponding to $f \in C^{\infty}_{\mathrm{o}}(S^{n-1})$ and
$\hat{F}$ is the cooriented subspace parallel to the facet $F$.
Using the identification of $\mathbb{P}_+^{\vee}(V)$ with
$\mathbb{P}_+(V^*)$, we can identify the line bundle $S$ over
$\mathbb{P}_+^{\vee}(V)$ with the line bundle $E \otimes
\mathscr{D}(V)$ over $\mathbb{P}_+(V^*)$, where $E$ is the line
bundle from part (a) of the proof. Thus, there is a canonical
isomorphism of Fr\'echet spaces
\begin{equation} \label{nun1inv2}
C^{\infty}_{\mathrm{o}}(\mathbb{P}_+(V^*),E) \otimes
\mathscr{D}(V) \rightarrow
C^{\infty}_{\mathrm{o}}(\mathbb{P}_+^{\vee}(V),S)
\end{equation}
which is $\mathrm{GL}(V)$ equivariant. Together, (\ref{nun1inv1})
and (\ref{nun1inv2}) determine a continuous $\mathrm{GL}(V)$
equivariant map
\begin{equation} \label{inviso2}
C^{\infty}_{\mathrm{o}}(\mathbb{P}_+(V^*),E) \otimes
\mathscr{D}(V) \rightarrow \mathbf{Val}_{n-1}^{\infty}
\end{equation}
whose image is dense by the Irreducibility Theorem and closed by
the Casselmann--Wallach theorem \textbf{\cite{cassel89}}. \hfill
$\blacksquare$

\vspace{0.3cm}

Now let $M$ again be a compact smooth manifold and $\pi: E \rightarrow M$ a line bundle over $M$.
In the same way the Poincar\'e duality map motivated the definition of generalized valuations, the pairing
(\ref{pairing}) motivates the following definition of the space of \emph{generalized sections} of $E$:
\[C^{-\infty}(M,E) = C(M,E^* \otimes \mathscr{D}M)^*.  \]
Note that the pairing (\ref{pairing}) yields a canonical embedding
\begin{equation} \label{canembsec}
C^{\infty}(M,E) \hookrightarrow C^{-\infty}(M,E).
\end{equation}

Finally, we are in a position to proof the main result of this appendix.

\vspace{0.2cm}

\noindent {\it Proof of Theorem \ref{genvalchar}.} The first part of the statement was already established in the proof
of Theorem \ref{thma1}. In order to prove the second statement, we first compute explicitly the Poincar\'e duality
map $\mathbf{Val}_1^{\infty}\! \rightarrow
\mathbf{Val}_1^{-\infty}$. To this end, let $\mu \in \mathbf{Val}_1^{\infty}$ and $\nu \in \mathbf{Val}_{n-1}^{\infty}$
be given by
\[\mu(K) = V(K,L,\ldots,L) \qquad \mbox{and} \qquad \nu(K) = V(K,\ldots,K,M)  \]
for some strictly convex bodies $L, M \in \mathcal{K}^n$ with smooth boundary. By the Irreducibility Theorem,
linear combinations of valuations of this form are dense in $\mathbf{Val}_1^{\infty}$ and $\mathbf{Val}_{n-1}^{\infty}$, respectively.
From (\ref{specialprod}), it follows that
\[<\!\mu,\nu\!> = \frac{1}{n} V(-L,\ldots,-L,M)\,V_n = \frac{V_n}{n^2}\int_{S^{n-1}} h(M,u)\,dS_{n-1}(-L,u).   \]
Thus, if $\bar{f} \in C^{\infty}_{\mathrm{o}}(\mathbb{P}_+(V^*),E^* \otimes
\mathscr{D}\,\mathbb{P}_+(V^*))$ is the section corresponding to the valuation $\mu \in \mathbf{Val}_1^{\infty}$ according to (\ref{inviso1})
and $\bar{g} \in C^{\infty}_{\mathrm{o}}(\mathbb{P}_+(V^*),E) \otimes \mathscr{D}(V)$ is the section
corresponding to $\nu \in \mathbf{Val}_{n-1}^{\infty}$ according to (\ref{inviso2}), then
\[\int_{S^{n-1}} h(M,u)\,dS_{n-1}(-L,u) = \int_{\mathbb{P}_+(V^*)} [\bar{f} \circ a,\bar{g}],  \]
where $a: \mathbb{P}(V^*) \rightarrow \mathbb{P}(V^*)$ denotes the antipodal involution on $\mathbb{P}(V^*)$, that is, the change of orientation. Now, noting that, by (\ref{inviso2}) and the fact that $\mathscr{D}(V)^* \otimes \mathscr{D}(V)$ is trivial,
\[\mathbf{Val}_1^{-\infty} = \left (\mathbf{Val}_{n-1}^{\infty} \right )^* \otimes \mathscr{D}(V) \cong C^{-\infty}_{\mathrm{o}}(\mathbb{P}_+(V^*),E^* \otimes
\mathscr{D}\,\mathbb{P}_+(V^*)) \]
and using again the isomorphism (\ref{inviso1}), we see that the Poincar\'e duality map induces a map
$C^{\infty}_{\mathrm{o}}(\mathbb{P}_+(V^*),E^* \otimes
\mathscr{D}\,\mathbb{P}_+(V^*)) \rightarrow C^{-\infty}_{\mathrm{o}}(\mathbb{P}_+(V^*),E^* \otimes
\mathscr{D}\,\mathbb{P}_+(V^*))$, given by
\[\bar{f} \mapsto \frac{1}{n^2} \bar{f} \circ a.  \]
Here we have used the embedding (\ref{canembsec}). This map obviously extends to an isomorphism of topological vector spaces
equipped with weak topologies
\[C^{-\infty}_{\mathrm{o}}(\mathbb{P}_+(V^*),E^* \otimes
\mathscr{D}\,\mathbb{P}_+(V^*)) \rightarrow C^{-\infty}_{\mathrm{o}}(\mathbb{P}_+(V^*),E^* \otimes
\mathscr{D}\,\mathbb{P}_+(V^*)) \cong \mathbf{Val}_1^{-\infty}.  \]
However, if we endow $V$ with a Euclidean structure, this map becomes an isomorphism
\[C^{-\infty}_{\mathrm{o}}(S^{n-1}) \rightarrow \mathbf{Val}_1^{-\infty}  \]
which, when restricted to smooth functions, is just given by
\[ f \mapsto \left (K \mapsto \frac{1}{n^2}\int_{S^{n-1}}f(-u)\,h(K,u)\,du \right).  \]
Clearly, this implies the desired statement.
\hfill $\blacksquare$

\pagebreak

\noindent {{\bf Acknowledgments} The work of A.\ Berg, L.\
Parapatits, and F.E.\ Schuster was supported by the European
Research Council (ERC), within the project ``Isoperimetric
Inequalities and Integral Geometry", Project number: 306445. L.\ Parapatits
was also supported by the ETH Zurich Postdoctoral Fellowship Program and the Marie Curie Actions
for People COFUND Program.

\begin{small}

\[ \begin{array}{ll} \mbox{Vienna University of Technology \phantom{wwwWWWW}} & \mbox{ETH Zurich} \\
\mbox{Institute of Discrete Math.\ and Geometry } &
\mbox{Department of Mathematics} \\
\mbox{Wiedner Hauptstrasse 8--10} & \mbox{R\"amistrasse 101} \\
\mbox{1040 Vienna, Austria} & \mbox{8092 Zurich, Switzerland}
\\[0.1cm]
\mbox{astrid.berg@tuwien.ac.at} & \mbox{lukas.parapatits@math.ethz.ch} \\
\mbox{franz.schuster@tuwien.ac.at} & \\
\mbox{manuel.weberndorfer@tuwien.ac.at} &
\end{array}\]

\end{small}

\end{document}